\documentclass{paper1}

\usepackage{ctrfont}
\usepackage{natbib}

\usepackage{graphicx}
\usepackage{subcaption}
\usepackage{epstopdf}

\usepackage{amsmath, rotating, amsfonts}

\usepackage{url}

\usepackage{color}
\usepackage{array}

\newcommand{\doublerightharpoonup}{%
  \rightharpoonup\mkern-10mu\rightharpoonup%
}
\newcolumntype{L}[1]{>{\raggedright\let\newline\\\arraybackslash\hspace{0pt}}m{#1}}
\newcolumntype{C}[1]{>{\centering\let\newline\\\arraybackslash\hspace{0pt}}m{#1}}
\newcolumntype{R}[1]{>{\raggedleft\let\newline\\\arraybackslash\hspace{0pt}}m{#1}}

\title{Adjoint-based gradient estimation from \\ gray-box solutions of unknown conservation laws} 

\author{Han CHEN\footnote[1]{Aerospace Computational Design Laboratory, Massachusetts Institute of Technology} \and Qiqi WANG\footnotemark[1]}

\begin{document}
\maketitle

\begin{abstract}
Many engineering applications can be formulated as optimizations constrained by conservation laws.
Such optimizations can be efficiently solved by the adjoint method, which computes the gradient of
the objective to the design variables.
Traditionally, the adjoint method has not been able to be
implemented in many ``gray-box'' conservation law 
simulators. In gray-box simulators, the analytical and numerical
form of the conservation law is unknown, but the full solution of relevant flow quantities is
available.
In this paper, we consider the case where the flux function is unknown.
This article introduces a method to estimate the gradient by
inferring the flux function from the solution, and then solving the adjoint equation of
the inferred conservation law.
This method is demonstrated in the sensitivity analysis of two flow problems.
\end{abstract}

\section{Introduction}
\label{background}
Optimization problems are of great interest in the engineering community. We consider an 
optimization problem to be
constrained by conservation laws.
For example, oil reservoir simulations may employ PDEs of various flow models, 
in which different fluid phases and components satisfy a set of conservation laws.
Such simulations can be used to facilitate the oil reservoir management,
including optimal well placement \cite{adjoint well placement} 
and optimal production control \cite{water flooding control,first reservoir opt}.
Another example is the cooling of turbine airfoils. 
We are interested in optimizing the interior flow path
of turbine airfoil cooling to minimize pressure loss
\cite{ubend rans opt 1, ubend rans opt 2}.\\

In many cases, such simulations can be computationally costly, potentially 
due to the complex computational models involved, and
large-scale time and space discretization. Furthermore, 
the dimensionality of the design space, $d$, can be high. 
For example, in oil reservoir simulations, the well pressure can be controlled at each well
individually, and they can vary in time.
To parameterize the well pressure, we require $Nm$ number of design variables, where $N$ is the number of wells,
and $m$ is the number of parameters, to describe the variation in time for each well.
Similarly, in turbine airfoil cooling,
the geometry of the internal flow path can also be parameterized by many variables.
Optimizing a high-dimensional design, $c$,
can be challenging.
A tool to enable efficient high-dimensional optimization is adjoint sensitivity analysis
\cite{adjoint}
, which efficiently computes the gradient of the objective to the design variables. 
The continuous adjoint method solves 
a continuous adjoint equation derived from the conservation law, which 
requires the PDE of the conservation law.
The discrete adjoint method solves a discrete adjoint equation derived 
from the numerical implementation of the conservation law,
which requires the simulator's numerical implementation.
Adjoint automatic differentiation applies the chain rule to every elementary arithmetic operation of the simulator,
which requires accessing and modifying the simulator's source code.
The adjoint methods has been applied to the sensitivity analysis of many 
problems, such as aerodynamic design \cite{adjoint aerodynamics} , 
oil reservoir history matching \cite{adjoint history matching} and optimal control 
\cite{adjoint reservoir optimal control, adjoint well placement}.
\\

We are mainly interested in gray-box simulations. By gray-box, we mean a conservation law simulation
 without the adjoint method implemented. Furthermore, we are not able to implement the adjoint method
when the governing PDE for the conservation law and its numerical implementation is unavailable:
for example, when the source code is proprietary or legacy.
Another defining property of gray-box simulation is that it can provide the space-time solution of the conservation law.
If the simulation solves for time-independent problems, a gray-box simulation should be able to 
provide the steady state solution.
In contrast, we define a simulator to be a blackbox, if neither the adjoint nor the solution is available.
The only output of such simulations is the value of the objective function to be optimized.
If the adjoint method is implemented or is able to be implemented,
we call such simulations open-box.
We summarize their differences in Table \ref{tab: boxes}.\\

\begin{center}
    \captionof{table}{Comparison of black-box, gray-box, and open-box simulations: the availability
                      of information}
    \label{tab: boxes}
    \begin{tabular}{ccccc}
        \hline
                   & PDE and implementation & {Space (or space-time) solution} & 
                   Adjoint\\[5pt] \hline
        Black-box  & No       & No   & No  \\ \hline
        Gray-box   & No
                   & Yes    & No   \\ \hline
        Open-box   & Yes    & Yes         &   Yes      \\ \hline
    \end{tabular}
\end{center}

We are also interested in the scenario where the space-time or spatial state variables
can be measured experimentally. For example, the Schlieren imaging technique is widely used to 
visualize transparent flows by the deflection of light, due to the refractive index gradient. 
The imaging can be used to reconstruct the distribution of the flow density 
\cite{schlieren reconstruct}. Another example is the plasma diagnostics which includes a large pool of
experimental techniques to measure the plasma properties \cite{francis chen}.
For example, the Langmuir probe and the magnetic probe can measure the local temperature,
the ion density, and the velocity of a plasma. 
The spectroscopic methods can measure line-integrated quantities such as the 
number density of certain ions species. Those techniques 
are widely used to monitor and control plasmas in real time. 
Such experimental techniques provides abundant data to reconstruct the state variables of the fluid
with a certain fidelity. However, similar to the case of gray-box simulation,
the experiments do not provide the sensitivity information directly.
This paper will focus on the gray-box simulation, and assume that the state variables are generated
from simulations. But our framework is also applicable
when the state variables are generated experimentally.\\

Depending on the type of simulation involved, we may choose different optimization methods.
If the simulation is black-box, we may use derivative-free optimization methods.
Derivative-free methods require only the availability of the objective function value, but
not the derivative information \cite{gradfreereview}.
Such methods are popular because they are easy to use. However, when the dimension of the design space 
increases, derivative-free methods may suffer from the curse of dimensionality.
The curse of dimensionality refers to problems caused by the rapid increase in the search 
volume associated
with adding an extra dimension in the search space. 
The resulting increase of search volume increases the number of objective evaluation required. 
It is not uncommon to encounter tens or hundreds of dimensions in real life engineering problems,
making derivative-free methods computationally expensive.\\

If the simulation is open-box, we may use gradient-based optimization methods.
Such methods use the gradient information to locate a local optimum.
A well-known example is the quasi-Newton method \cite{quasiNewton}. 
When the design space is high-dimensional, gradient-based methods can require fewer
objective evaluations to converge than derivative-free methods.
Let $c$ be the design variables, $u$ be the flow solution, and $J$ be the objective value.
Gradient-based optimization methods require $\frac{d J}{d c} = \frac{\partial J}{\partial c} + \frac{\partial J}{\partial u}
\frac{d u}{d c}$, the total derivative of the objective to the design variables. 
The term $\frac{dJ}{dc}$ can be evaluated efficiently by using adjoint methods.
The continuous adjoint method derives the continuous adjoint equation from the PDE of the simulation through the
method of the Lagrange multiplier; therefore, it requires knowing the PDE of the simulator. 
The discrete adjoint method applies variational analysis to the discretized PDE, and it requires knowing the
discretized PDE of the simulator.
The discrete adjoint method can be implemented by using automatic differentiation (AD).
AD decomposes the simulation into elementary arithmetic operations and functions, and
then applies the chain
rule to compute the gradient.
Because adjoint methods require access to the PDE and its discretization, they cannot be directly
applied to gray-box simulations.\\

If the simulation is gray-box, we cannot apply the adjoint method to the gray-box conservation law
to compute the gradient. In current practice, gray-box simulations are often treated as black-box simulations, to which 
adjoint-based optimization is not applicable.
The gray-box simulation is viewed as a calculator for the objective, while its space-time solution is neglected.
This paper introduces a method to enable adjoint-based optimization on gray-box simulations. We propose
a two-step procedure to estimate the objective's gradient by using the gray-box simulation's space-time solution.
In the first step, we infer the conservation law governing the simulation by using its space-time solution. In
the second step, we apply the adjoint method to the inferred conservation law to estimate the gradient.\\

In the remainder of the paper, Section \ref{prototypes} defines the particular type of
gray-box models considered in this paper. Section \ref{infer} explains why it can be feasible to
infer the governing PDE from the output of the gray-box model, i.e. the space-time solution.
Section \ref{inverse} describes how the PDE can be inferred. 
Section \ref{numerical example} demonstrates that our method works on a 1-D porous media flow problem and
a 2-D Navier-Stokes flow problem.

\section {Optimization constrained by gray-box conservation laws: \\Problem definition}
\label{prototypes}

We consider the case in which the governing PDE is a conservation law.
For example, time-dependent conservation laws can be written as
\begin{equation}
    \dot{\boldsymbol{u}} + \nabla \cdot \overset{\doublerightharpoonup}{\boldsymbol{F}}
    (\boldsymbol{u})
    = \boldsymbol{q}(\boldsymbol{u},c)
    \label{first equation}
\end{equation}
where $\boldsymbol{u}$ is a vector representing flow quantities,
$\dot{\boldsymbol{u}}$ is the derivative of $\boldsymbol{u}$ with respect to time,
$c$ represents the design variables,
$t\in[0,T]$ is time;
and $x\in \Omega \subseteq \mathbb{R}^{n}$ is the spatial coordinate.
$\Omega$ may depend on the design variables, $c$.
$\overset{\doublerightharpoonup}{\boldsymbol{F}}$ is the flux tensor.
$\boldsymbol{q}$ is a source vector that may also depend on $c$.
The boundary and initial conditions are known.
The discretized space-time solution of Eqn.\eqref{first equation} given by a gray-box simulation
is written as $\hat{u}(t_i, \mathbf{x}_i; c)\,, i=1,\cdots,N$, where
$t=\left\{t_1,\cdots, t_N\right\}$ indicates the discretized time, and 
$\mathbf{x}_i$ indicates the spatial discretization at time $t_i$.\\

If we are only interested in the steady state solution of the conservation law, the flow solution satisfies
\begin{equation}
    \nabla \cdot 
    \overset{\doublerightharpoonup}{\boldsymbol{F}}(\boldsymbol{u}_{\infty}) 
    = \boldsymbol{q}_i(\boldsymbol{u}_{\infty},c)
    \label{first equation steady}
\end{equation}
where $\boldsymbol{u}_{\infty}$ is the converged solution of 
of Eqn.\eqref{first equation} at $t\rightarrow \infty$.\\

In many cases, Eqns.\eqref{first equation} or \eqref{first equation steady} are unknown. 
In this paper,
we consider the case where the flux $\overset{\doublerightharpoonup}{\boldsymbol{F}}$ is unknown. 
We provide two examples to illustrate conservation laws with unknown fluxes. 
The first example is a 1-D two-phase flow in porous media \cite{Buckley Leverett}. 
The governing equation can be written as
\begin{equation}
    \frac{\partial u}{\partial t} + \frac{\partial }{\partial x} F(u) = c
\end{equation}
where $x\in[0,1]$ is the space domain;
$u = u(t,x)$, $0\le u\le 1$, is the saturation of phase I (e.g. water), 
and $1-u$ is the saturation of phase II (e.g. oil);
$c=c(t,x)$ is the design variable. $c>0$ models the injection of phase I replacing phase II; and $c<0$
represents the opposite.
$F$ is an unknown flux that depends on the properties of the porous media and the fluids.\\

The second example is a 2-D Navier-Stokes flow for a fluid with an unknown state equation. 
Let $\rho$, $u$, $v$, $E$, and $p$ denote the density, Cartesian velocity components, 
total energy, and pressure.
The governing equation is
\begin{equation}
    \frac{\partial}{\partial t}
    \begin{pmatrix}
        \rho \\ \rho u \\ \rho v\\ \rho E
    \end{pmatrix}
    + \frac{\partial}{\partial x} 
    \begin{pmatrix}
        \rho u\\
        \rho u^2 + p - \sigma_{xx}\\
        \rho uv - \sigma_{xy}\\
        u(E\rho+p) - \sigma_{xx} u - \sigma_{xy} v
    \end{pmatrix}
    + \frac{\partial}{\partial y}
    \begin{pmatrix}
        \rho v\\
        \rho uv-\sigma_{xy}\\
        \rho v^2+p-\sigma_{yy}\\
        v(E\rho+p) - \sigma_{xy} u -\sigma_{yy}v
    \end{pmatrix} 
    = \boldsymbol{0}
    \label{NSeqn}
\end{equation}
where
\begin{equation}\begin{split}
    \sigma_{xx} &= \mu \left(2 \frac{\partial u}{\partial x} - \frac{2}{3} \left(\frac{\partial u}{\partial x} 
    + \frac{\partial v}{\partial y}\right)\right)\\
    \sigma_{yy} &= \mu \left(2 \frac{\partial v}{\partial y} - \frac{2}{3} \left(\frac{\partial u}{\partial x} 
    + \frac{\partial v}{\partial y}\right)\right)\\
    \sigma_{xy}&=\mu\left(\frac{\partial u}{\partial y} + \frac{\partial v}{\partial x}\right)
\end{split}\end{equation}
The pressure is governed by an unknown state equation
\begin{equation}
    p = p(U, \rho)
    \label{state equation}
\end{equation}
where $U$ denotes the internal energy per volume,
\begin{equation}
    U = \rho\left(E-\frac{1}{2}(u^2+v^2)\right)\,.
\end{equation}
The state equation depends on the property of the fluid.\\

For time-dependent problems, we are interested in solving the optimization problem
\begin{equation}
    \min_{c}J = \min_{c}\int_0^T \int_\Omega j(\boldsymbol{u},c) \textrm{d}\mathbf{x}\textrm{d}t
    \label{obj prototype}
\end{equation}
where $\boldsymbol{u}$ satisfies Eqn.\eqref{first equation}
whose flux, $\overset{\doublerightharpoonup}{\boldsymbol{F}}$, is unknown.
Notice that $j$ can depend explicitly on c, while $u$ depends implicitly on $c$.
If we are interested in the steady state, we optimize
\begin{equation}
    \min_{c}J =\min_{c} \int_\Omega j(\boldsymbol{u}_{\infty},c) \textrm{d}\mathbf{x}
\end{equation}
where $\boldsymbol{u}_{\infty}$ satisfies Eqn.\eqref{first equation steady}
whose flux, $\overset{\doublerightharpoonup}{\boldsymbol{F}}$, is unknown.

\section{Infer conservation laws by the space-time solution}
\label{infer}
Because Eqns. \eqref{first equation} or \eqref{first equation steady}
haves an unknown flux, the adjoint method cannot be applied to evaluate 
$\frac{dJ}{dc}$. To enable the adjoint method in such a scenario, we will infer
the unknown flux.\\

We use the space-time solution of the gray-box simulations 
to infer the flux.
There are several benefits by using the space-time solution \cite{hanmaster, ecmor}.
Firstly, in conservation law simulations, the flow quantities only depend on the flow
quantities in a previous time inside a domain of dependence.
When the timestep is small,
the domain of dependence can be small, as well. For example, for scalar conservation laws 
without exogenous control,
we can view solving the conservation law for one timestep $\Delta t$ at a spatial location
as a mapping 
$\mathbb{R}^{\omega_{\Delta t}} \rightarrow \mathbb{R}$, where $\omega_{\Delta t}\in \Omega$ 
is the domain of dependence. 
By applying such mapping repeatedly to all $x\in \Omega$ and $t\in[0,T]$
(in addition to the boundary and initial conditions), we perform 
a space-time simulation of the conservation law. 
Generally, the size of $\omega_{\Delta t}$ is small.
Therefore, in the discretized simulation of conservation laws,
the number of discretized flow variables involved in $\omega_{\Delta t}$ is small, as well,
making it feasible to infer the mapping.
\\

\begin{figure}\begin{center}
    \includegraphics[width=0.5\textwidth, natwidth=100,natheight=100]{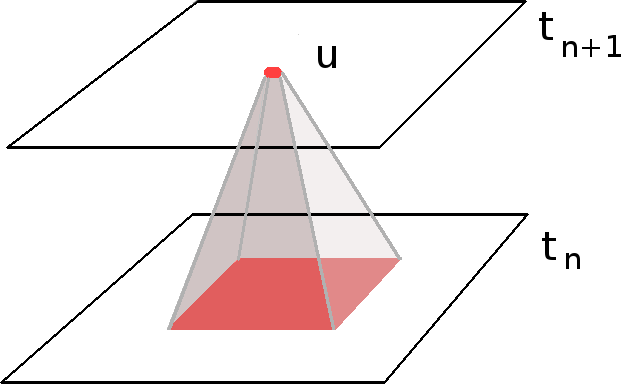}
    \caption{Domain of dependence: in conservation law simulations,
             the flow quantities at a given location
             depend on the flow quantities at an older time only within its
             domain of dependence. The two planes in this figure indicates the spatial 
             solution at two adjacent timesteps.
             The domain of dependence can be much smaller
             than the overall spatial domain when the timestep is reasonably small.}
    \label{locality}
\end{center}\end{figure}

Secondly, the space-time solution at almost {every} space-time grid point can be viewed
as a sample for the mapping $\mathbb{R}^\omega \rightarrow \mathbb{R}$. 
Because the number of space-time grid points in gray-box simulations is generally large,
we have a large number of samples to infer the mapping. 
In Eqn.\eqref{first equation} and \eqref{first equation steady},
such a mapping is determined by the flux.
Therefore, we will have a large number of samples to infer the flux,
making the inference potentially accurate.\\

Thirdly, in many optimization problems, the design space is high dimensional only because
the design is space- and/or time-dependent. In order to parameterize the space-time dependent
design, a large number of design variables will be employed. However, 
the flow quantities only depend on the design variables in the domain of dependence.
Therefore, even if the overall number of design variables is high, the number of design variables
involved in the mapping is limited, making the inference problem potentially
immune to the design space dimensionality.\\

Therefore, we propose to infer the flux in Eqn.\eqref{first equation} or
\eqref{first equation steady}. The inferred
PDE should yield a solution, $\tilde{\boldsymbol{u}}$, that matches 
the solution of the gray-box simulation, $\boldsymbol{u}$.
A simulator of the inferred PDE is called the {twin model}.

\section{Twin model inference as an optimization problem}
\label{inverse}
Conventionally, we have a given PDE, and want to compute
its space-time solution. However, in a twin model, we want to infer the PDE
to match a given space-time solution.
Finding a suitable PDE, specifically a suitable flux function, can be viewed
as an inverse problem, which can be solved by optimization.
We define a metric for the mismatch of the space-time solutions.
Given the same inputs (design variables, initial conditions, and boundary conditions), 
a twin model should yield a space-time solution $\tilde{\boldsymbol{u}}$ such that 
$\tilde{\boldsymbol{u}}$ is close to $\boldsymbol{u}$. 
Suppose that the twin model and the 
primal model use the same discretization; we use the following expression to quantify
the mismatch:
\begin{equation}
    \mathcal{M} = \frac{1}{T}
    \sum_{i=1}^{N}\sum_{k=1}^{M} \left(\tilde{\boldsymbol{u}}_{ik} 
    - \boldsymbol{u}_{ik}\right)^2 \Delta t_k
    \left| \Delta \mathbf{x}_i \right|
    \label{minimizer twin model discrete}
\end{equation}
where $\left| \Delta \mathbf{x}_i \right|$ indicates the lengths (1-D), areas (2-D), or volumes (3-D) 
of the grid.
If the space- and/or time-grids are different, then a mapping $P$ from $u$ to $\tilde{u}$ is required.
In this case, the mismatch is defined by
\begin{equation}
    \mathcal{M} = \frac{1}{T}
    \sum_{i=1}^{N}\sum_{k=1}^{M} \left(\tilde{\boldsymbol{u}}_{ik} - P(\boldsymbol{u})_{ik}\right)^2 
    \Delta t_k
    \left| \Delta \mathbf{x}_i \right|
    \label{minimizer twin model discrete mapping}
\end{equation}
For the present examples, we assume that the grids are the same for the purpose of simplicity.\\

We parameterize the flux function and infer the parameterization that minimizes $\mathcal{M}$.
Let the parameterized flux function be $G(\tilde{\boldsymbol{u}}, \xi)$,
where $\xi$ values are the parameters.
The inference problem is stated as follows:\\

Solve
\begin{equation}
    \xi^* = 
    \arg\min_{\xi} \left\{
    \mathcal{M}
    + \lambda \|\xi\|^p  \right\}
    \label{objective twin model}
\end{equation}
where $\tilde{\boldsymbol{u}}$
is the discretized space-time solution of
\begin{equation}
    \dot{ \tilde{\boldsymbol{u}}} + \nabla \cdot
    G(\tilde{\boldsymbol{u}}, \xi)
    = \boldsymbol{q}(\tilde{\boldsymbol{u}},c)
    \label{first equation 2}
\end{equation}
 $\lambda\|\xi\|^p$
is an $L_p$ norm regularization, and $\lambda>0$.
The space-time discretization of Eqn.\eqref{first equation 2} is the same as the gray-box simulator.
\\


Similarly, if we are interested in the steady state solution,
we have
\begin{equation}
    \mathcal{M} = \frac{1}{T}
    \sum_{i=1}^{N} \left(\tilde{\boldsymbol{u}}_{i} - \boldsymbol{u}_{i}\right)^2
    \left| \Delta \mathbf{x}_i \right|
    \label{minimizer twin model discrete steady}
\end{equation}

The inference problem is stated as follows:\\

Solve
\begin{equation}
    \xi^* = 
    \arg\min_{\xi} \left\{
    \mathcal{M} 
    + \lambda \|\xi\|^p  \right\}
    \label{objective twin model steady}
\end{equation}
where $u$ is the discretized spatial solution of the gray-box simulation, and $\tilde{u}$
is the discretized spatial solution of
\begin{equation}
    \nabla \cdot
    G(\tilde{\boldsymbol{u}}, \xi)
    = \boldsymbol{q}(\tilde{\boldsymbol{u}},c)
    \label{first equation 2 steady}
\end{equation}

The match of space-time solutions does not ensure the match of flux functions.
The problem can be ill-posed to infer $F$ on a domain not covered by the gray-box solution $u$.
In other words, the basis for modeling the flux may be over-complete.
We will call a domain of $u$ ``excited'' if inferring $F$ is well-posed on that domain.
The ill-posedness
can be allieviated by basis selection. It has been shown that basis selection can be performed
by Lasso regularization corresponding to $p=1$ \cite{Lasso variable selection}.
We will set $p=1$ in this article. \\

Because the twin model is an open-box system, Eqns.\eqref{objective twin model} 
and \eqref{objective twin model steady} can be solved by adjoint gradient-based methods.
The parameterization of the twin model flux $G$ will be problem-dependent. We will discuss
this topic in Section \ref{numerical example}.

\section{Numerical examples}
\label{numerical example}
In this section, we demonstrate the twin model with two numerical examples.
\subsection{Gradient estimation for a 1-D porous media flow}
Consider a 1-D PDE
\begin{equation}
    \frac{\partial u}{\partial t} + \frac{\partial F(u)}{\partial x} = c\,\quad x\in[0,1]\; t\in[0,1]
    \label{BL eqn}
\end{equation}
with periodic boundary condition
\begin{equation}
    u(x=0) = u(x=1)
\end{equation}
and initial condition
\begin{equation}
    u(t=0) = u_0
\end{equation}
Eqn.\eqref{BL eqn} can be used to model 1-D, two-phase, porous media flow, where $u$
denotes the saturation of one phase. $0\le u\le 1$. $c=c(t,x)$ is a space-time dependent exogenous control.
The flux function $F(u)$ depends on the properties of the porous media and the fluids.
For example, the Buckley-Leverett equation models the flow driven by 
capillary pressure and Darcy's law \cite{Buckley Leverett}, whose flux is
\begin{equation}
    F(u) = \frac{u^2}{1+A(1-u)^2}
    \label{BL flux}
\end{equation}
where $A$ is a constant. In the following we assume the graybox simulation
solves Eqn.\eqref{BL eqn} and \eqref{BL flux} with $A=2$.
\\

Assuming that $F(u)$ is unknown, we will fit a twin model 
using the graybox simulation.
We parameterize the twin model according to Eqn.\eqref{first equation 2}
\begin{equation}
    \frac{\partial \tilde{u}}{\partial t} + \frac{\partial}{\partial x}\,
    \left(\sum_{k=1}^m \xi_k g_{k}(\tilde{u})\right) = c
    \label{twin model 3}
\end{equation}
with the same initial condition, boundary conditions, and exogenous control. 

Generally, the flux $F$ is a monotonic increasing function. To respect this fact, we
choose $g$ values as a family of sigmoid functions 
\begin{equation}
    g_k(\tilde{u}) = \left(\tanh\left(\frac{\tilde{u}- \eta_k}{\sigma}\right)+1\right) \left.\right/2
    \label{sigmoids}
\end{equation}
where $\eta_k$ and $\sigma$ are constants. Therefore, we just need to set $\xi\ge 0 $ 
to enforce the monotonicity of the flux.
We will use  a second-order finite volume discretization and 
Crank-Nicolson time integration scheme in the twin model.
\begin{figure}\begin{center}
    \includegraphics[width=7cm]{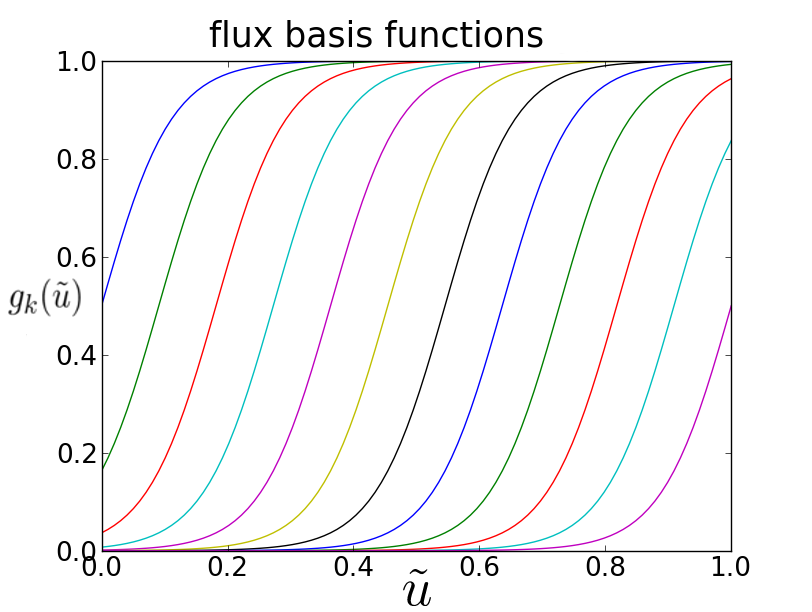}
    \caption{An example of the flux basis functions $g_k(\tilde{u})$. 
    These basis functions are distinguished by different colors.}
\end{center}\end{figure}

To infer the coefficients $\xi$, we minimize
\begin{equation}
    \int_0^1 \int_0^1 (u-\tilde{u})^2 dx dt+ \lambda \sum_{k=1}^m |\xi_k|
\end{equation}
We use the {low-memory Broyden-Fletcher-Goldfarb-Shannon} (BFGS) algorithm
\cite{LBFGS} for the minimization. 
L-BFGS approximates the Hessian using only the gradients at newer previous iterations,
and inverses the approximated Hessian efficiently using the Sherman-Morrison formula.
The gradient, $\frac{dJ}{d\xi_k}\,, k=1,\cdots, m$, 
is computed by an automatic differentiation module \textit{numpad} 
$[$Q. Wang, https://github.com/qiqi/numpad.git.$]$\\

\begin{figure}\begin{center}
    \includegraphics[height=8cm]{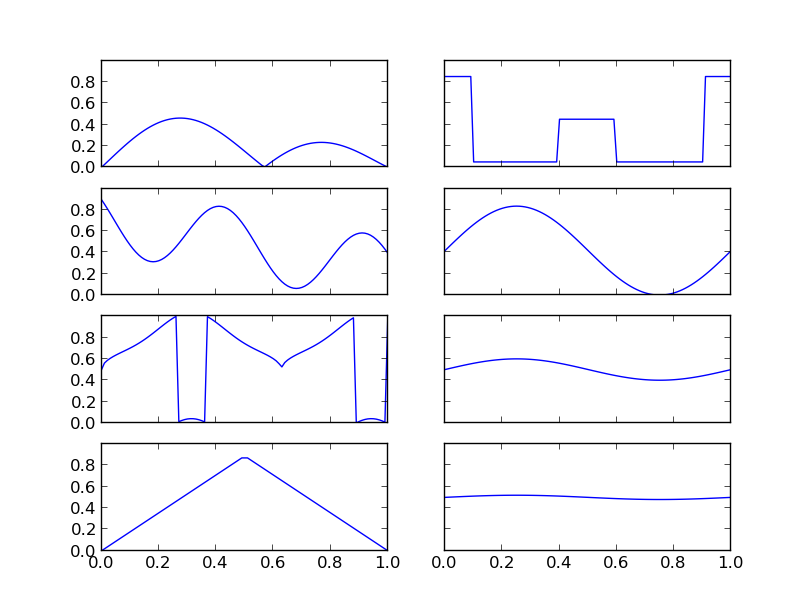}
    \caption{Initial condition $u_0(x)$. 
    We choose a diverse set of initial conditions to test the twin model method.}
    \label{fig:initial condition}
\end{center}\end{figure}

\begin{figure}
    \begin{subfigure}[b]{.95\textwidth}
          \centering
          \begin{subfigure}[b]{0.32\textwidth}
                  \centering
                  \includegraphics[width=4cm,height=3.2cm]{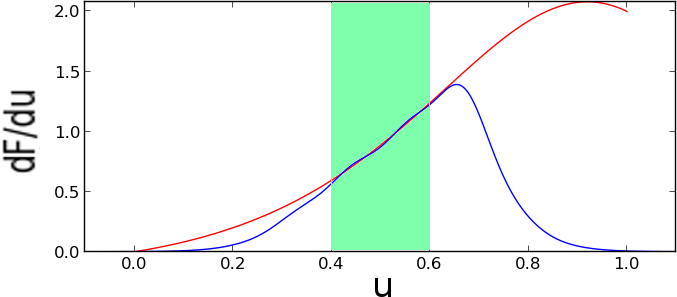}
          \end{subfigure}%
          \begin{subfigure}[b]{0.32\textwidth}
                  \centering
                  \includegraphics[width=4.3cm]{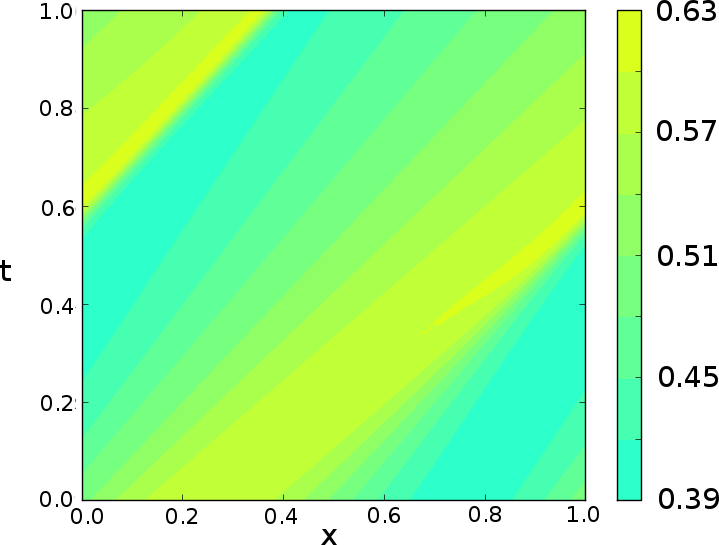}
          \end{subfigure}
          \begin{subfigure}[b]{0.32\textwidth}
                  \centering
                  \includegraphics[width=4.3cm]{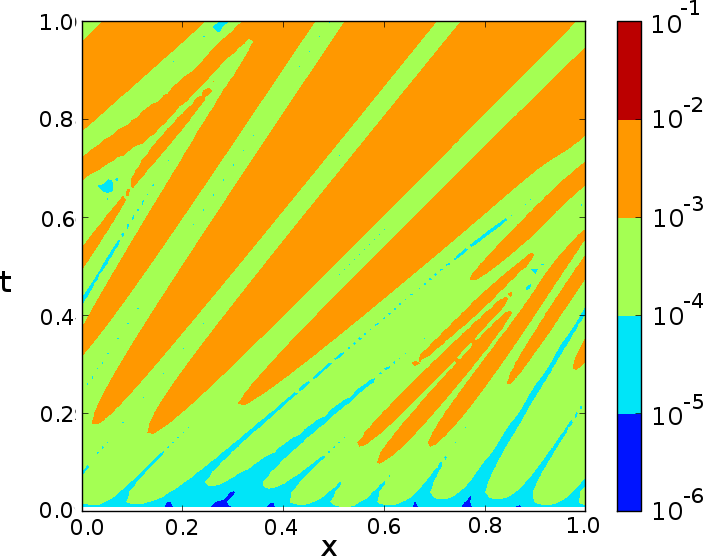}
          \end{subfigure}%
    \end{subfigure}
    \begin{subfigure}[b]{.95\textwidth}
          \centering
          \begin{subfigure}[b]{0.32\textwidth}
                  \centering
                  \includegraphics[width=4cm,height=3.2cm]{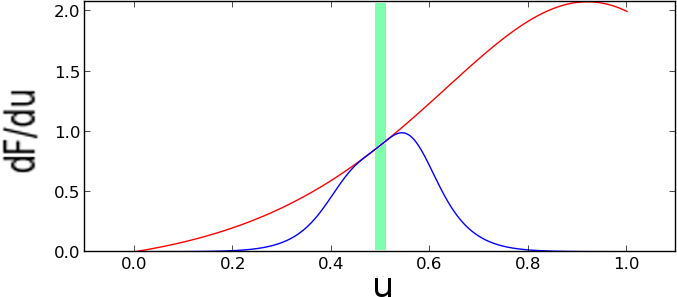}
          \end{subfigure}%
          \begin{subfigure}[b]{0.32\textwidth}
                  \centering
                  \includegraphics[width=4.3cm]{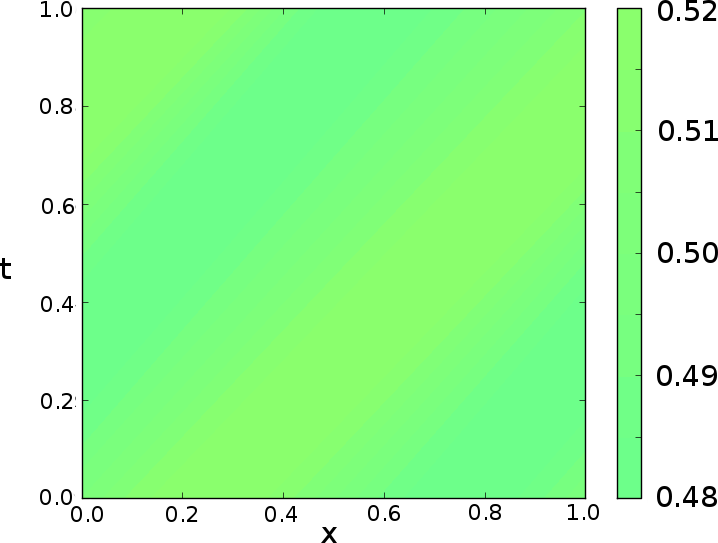}
          \end{subfigure}
          \begin{subfigure}[b]{0.32\textwidth}
                  \centering
                  \includegraphics[width=4.3cm]{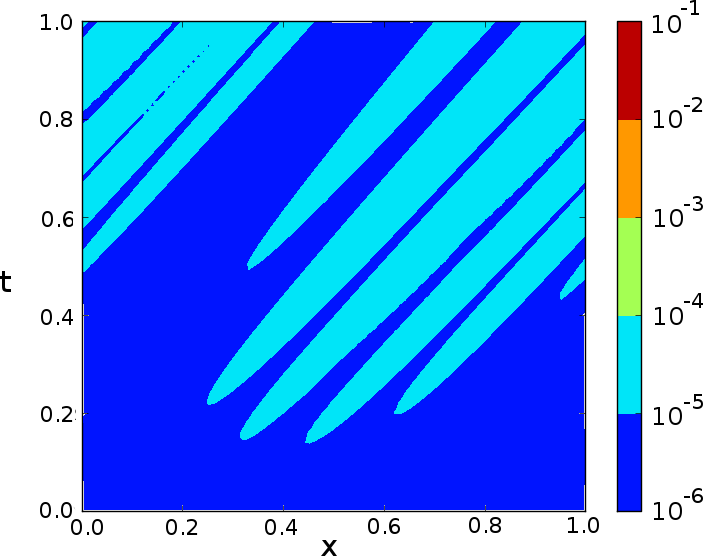}
          \end{subfigure}
    \end{subfigure}

    \caption{
             \label{fig:sol compare}
             The left column shows the gradient of $F$ (red line) with the gradient of $\tilde{F}$ 
             (blue line).  The middle column shows the space-time solutions of the gray-box model.
             The right column shows the solution mismatch $\left|u-\tilde{u}\right|$. 
             Each row corresponds to an initial condition.
             The excited domain, $[u_{\min}, u_{\max}]$, is indicated by the green region in
             the left column.
             When the excited domain is large, the twin model's flux 
             approximates the gray-box model's flux in a large range of $u$, but
             the solution mismatch tends to be less accurate in the exicted domain. 
             When the excited domain is small,
             the twin model's flux approximates the gray-box model's flux in a small range of $u$, 
             but the flux approximation tends to be more accurate in the excited domain.}
\end{figure}

We set $c=0$, and test the quality of the twin model given several initial conditions.
These initial conditions are shown in Fig. \ref{fig:initial condition}.
For different initial conditions, the excited domain will be different.
Let $u_{\max} = \max_{t,x}u(t,x)$ and $u_{\min} = \min_{t,x}u(t,x)$,
we expect the inferred flux to match the true flux within $[u_{\min},u_{\max}]$.
In Eqn.\eqref{first equation},
we can add a constant to the flux while yielding the same solution $u$. 
Therefore, we should compare $\frac{dF}{du}$ with $\frac{d\tilde{F}}{du}$, instead of comparing
$F$ with $\tilde{F}$. 
Fig.\ref{fig:sol compare} shows $\frac{dF}{du}$, $\frac{d\tilde{F}}{du}$, the space-time solution 
of the gray-box simulation, and the solution mismatch.\\

\begin{figure}
    \centering
    Flux gradient\hspace{2cm} Gray-box solution \hspace{2cm} Solution mismatch
    \begin{subfigure}[b]{.95\textwidth}
          \centering
          \begin{subfigure}[b]{0.32\textwidth}
                  \centering
                  \includegraphics[width=4cm,height=3.2cm]{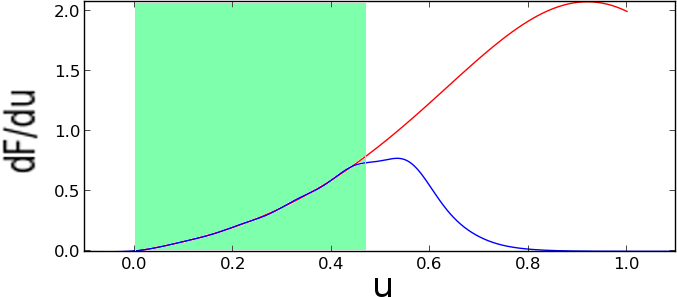}
          \end{subfigure}%
          \begin{subfigure}[b]{0.32\textwidth}
                  \centering
                  \includegraphics[width=4.3cm]{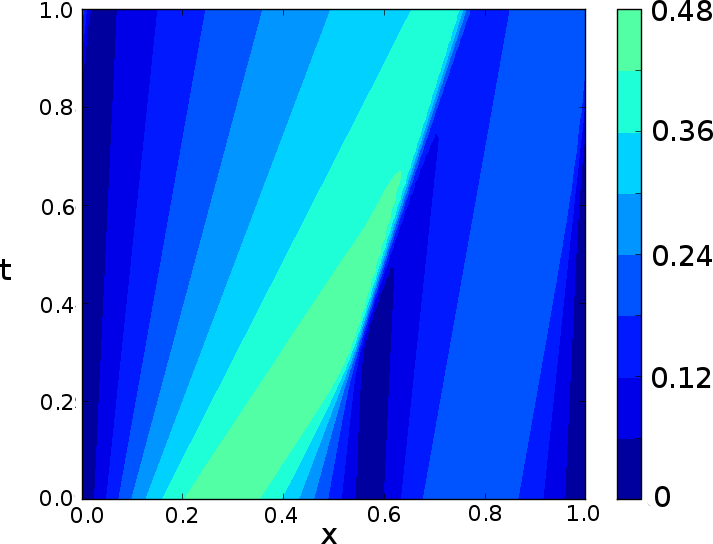}
          \end{subfigure}
          \begin{subfigure}[b]{0.32\textwidth}
                  \centering
                  \includegraphics[width=4.3cm]{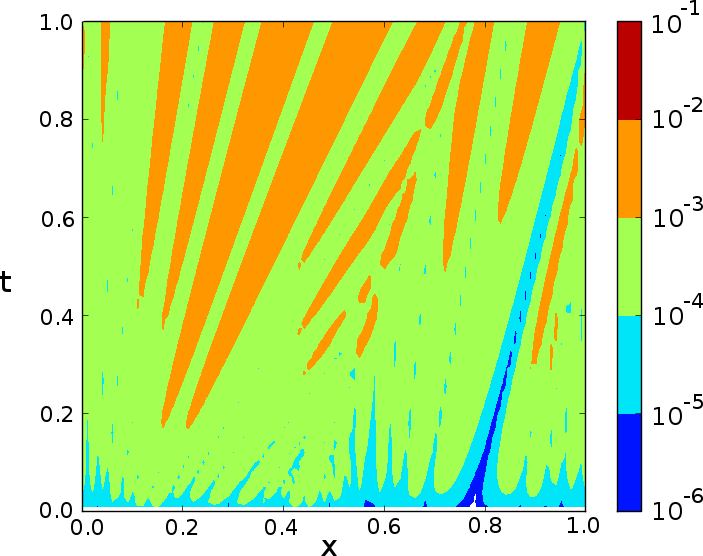}
          \end{subfigure}%
    \end{subfigure}
    \begin{subfigure}[b]{.95\textwidth}
          \centering
          \begin{subfigure}[b]{0.32\textwidth}
                  \centering
                  \includegraphics[width=4cm,height=3.2cm]{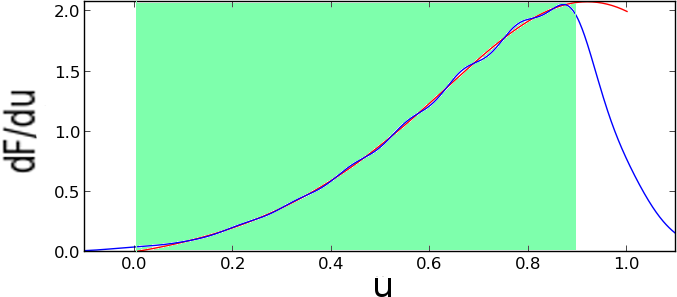}
          \end{subfigure}%
          \begin{subfigure}[b]{0.32\textwidth}
                  \centering
                  \includegraphics[width=4.3cm]{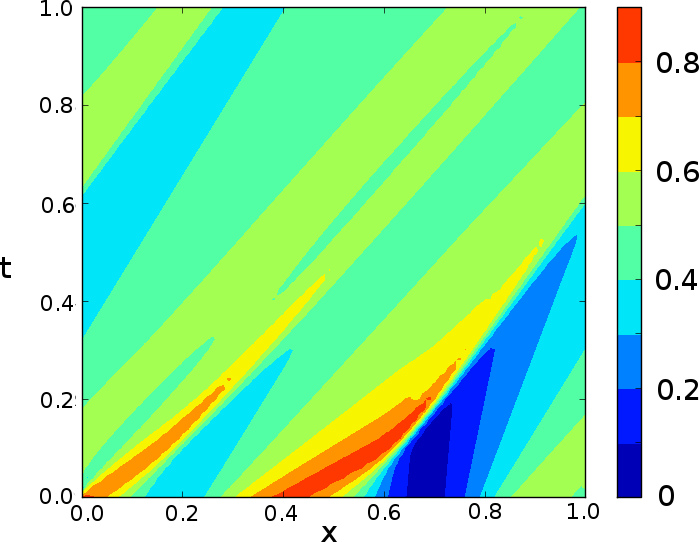}
          \end{subfigure}
          \begin{subfigure}[b]{0.32\textwidth}
                  \centering
                  \includegraphics[width=4.3cm]{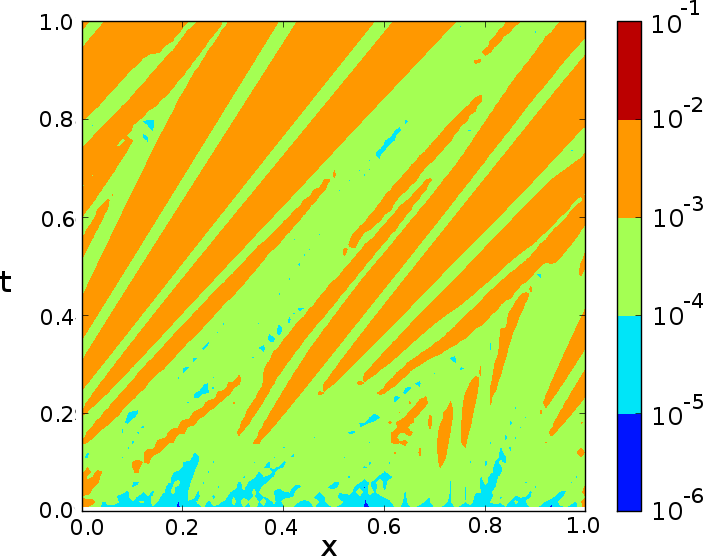}
          \end{subfigure}
    \end{subfigure}
    \begin{subfigure}[b]{.95\textwidth}
          \centering
          \begin{subfigure}[b]{0.32\textwidth}
                  \centering
                  \includegraphics[width=4cm,height=3.2cm]{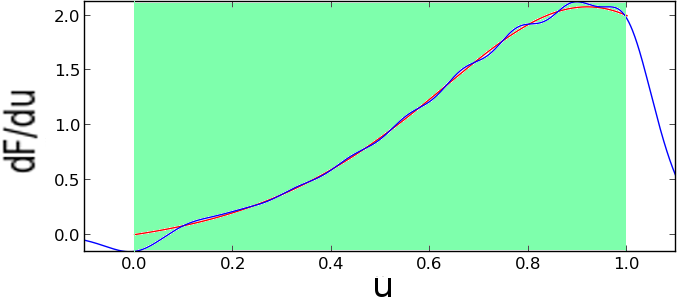}
          \end{subfigure}%
          \begin{subfigure}[b]{0.32\textwidth}
                  \centering
                  \includegraphics[width=4.3cm]{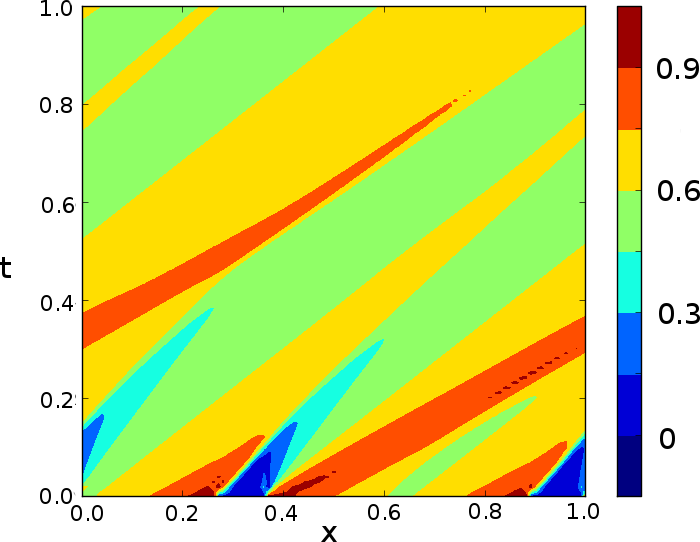}
          \end{subfigure}
          \begin{subfigure}[b]{0.32\textwidth}
                  \centering
                  \includegraphics[width=4.3cm]{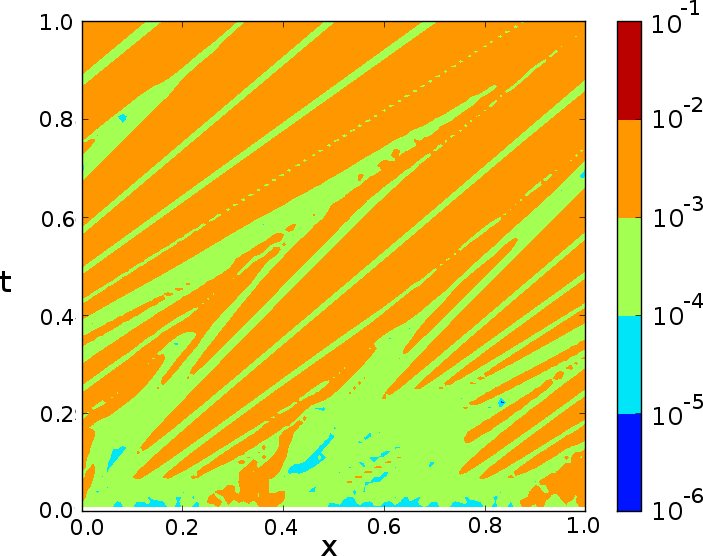}
          \end{subfigure}%
    \end{subfigure}
    \begin{subfigure}[b]{.95\textwidth}
          \centering
          \begin{subfigure}[b]{0.32\textwidth}
                  \centering
                  \includegraphics[width=4cm,height=3.2cm]{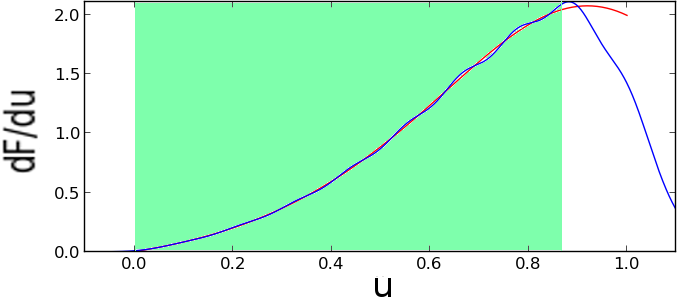}
          \end{subfigure}%
          \begin{subfigure}[b]{0.32\textwidth}
                  \centering
                  \includegraphics[width=4.3cm]{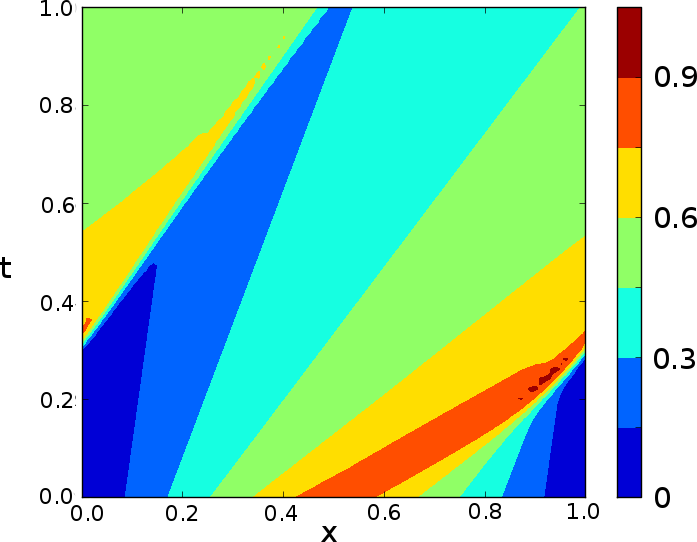}
          \end{subfigure}
          \begin{subfigure}[b]{0.32\textwidth}
                  \centering
                  \includegraphics[width=4.3cm]{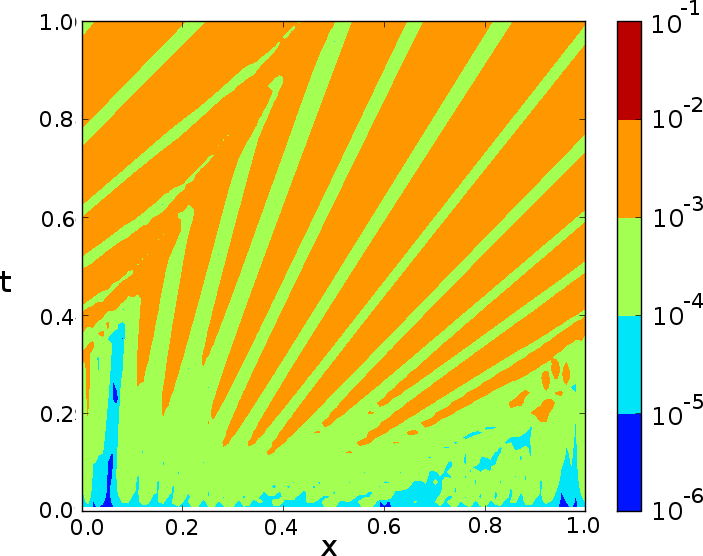}
          \end{subfigure}
    \end{subfigure}
    \begin{subfigure}[b]{.95\textwidth}
          \centering
          \begin{subfigure}[b]{0.32\textwidth}
                  \centering
                  \includegraphics[width=4cm,height=3.2cm]{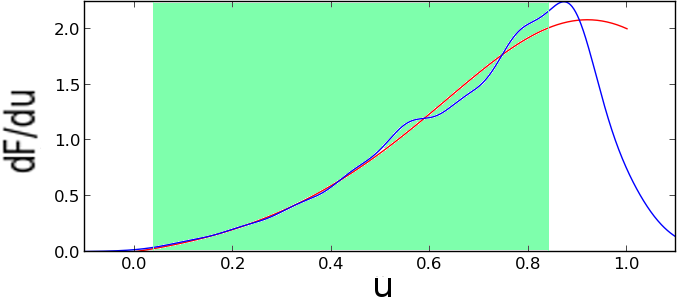}
          \end{subfigure}%
          \begin{subfigure}[b]{0.32\textwidth}
                  \centering
                  \includegraphics[width=4.3cm]{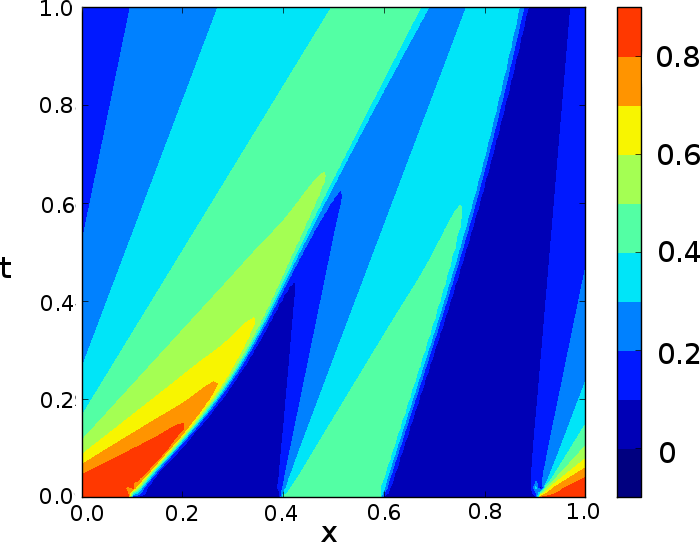}
          \end{subfigure}
          \begin{subfigure}[b]{0.32\textwidth}
                  \centering
                  \includegraphics[width=4.3cm]{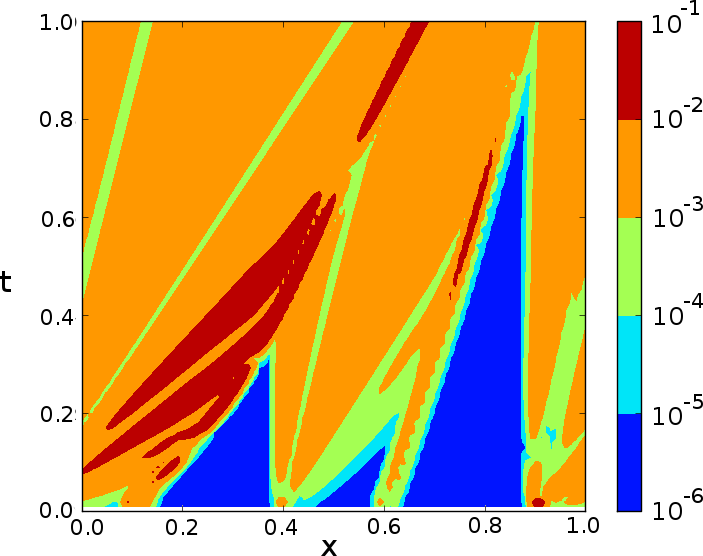}
          \end{subfigure}%
    \end{subfigure}
    \begin{subfigure}[b]{.95\textwidth}
          \centering
          \begin{subfigure}[b]{0.32\textwidth}
                  \centering
                  \includegraphics[width=4cm,height=3.2cm]{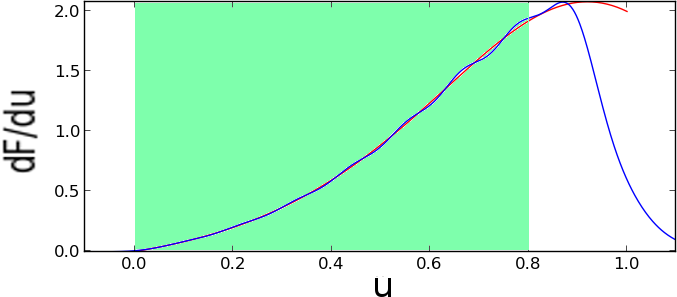}
          \end{subfigure}%
          \begin{subfigure}[b]{0.32\textwidth}
                  \centering
                  \includegraphics[width=4.3cm]{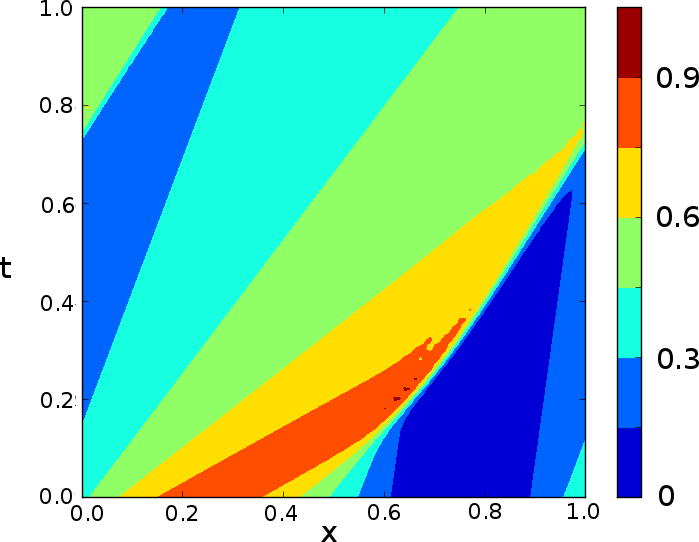}
          \end{subfigure}
          \begin{subfigure}[b]{0.32\textwidth}
                  \centering
                  \includegraphics[width=4.3cm]{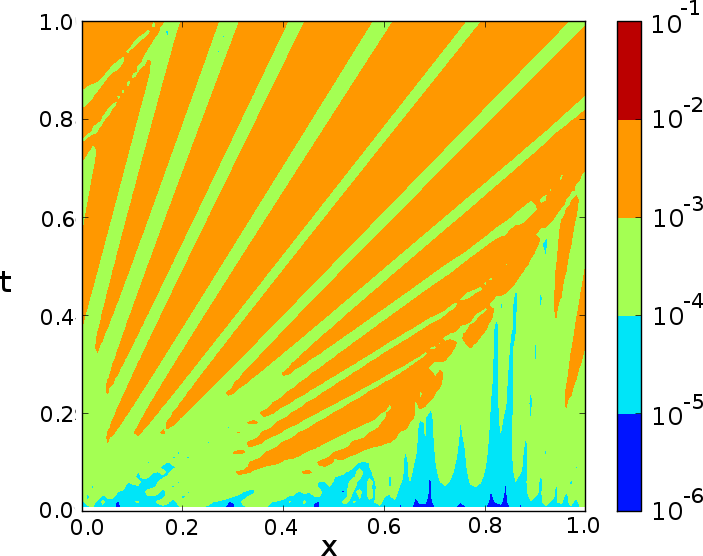}
          \end{subfigure}
    \end{subfigure}
\end{figure}

\begin{figure}
    \begin{subfigure}[b]{.95\textwidth}
          \centering
          \begin{subfigure}[b]{0.32\textwidth}
                  \centering
                  \includegraphics[width=4.3cm]{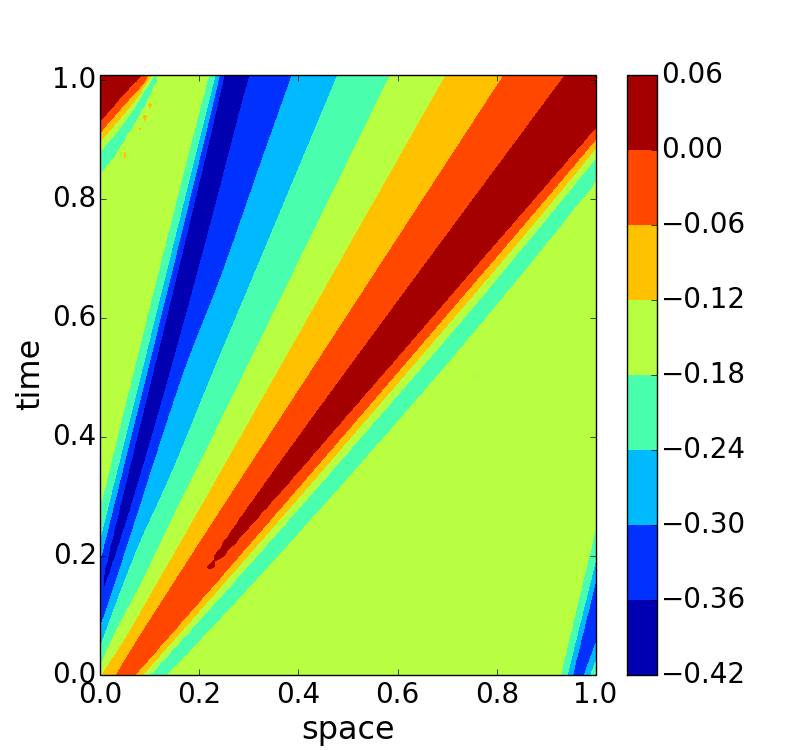}
          \end{subfigure}%
          \begin{subfigure}[b]{0.32\textwidth}
                  \centering
                  \includegraphics[width=4.3cm]{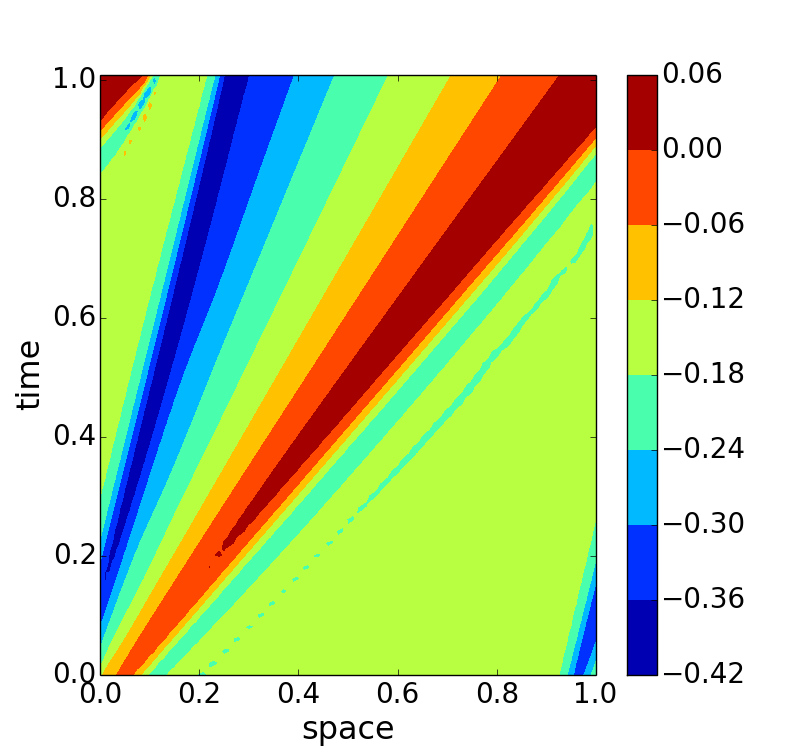}
          \end{subfigure}
          \begin{subfigure}[b]{0.32\textwidth}
                  \centering
                  \includegraphics[width=4.3cm]{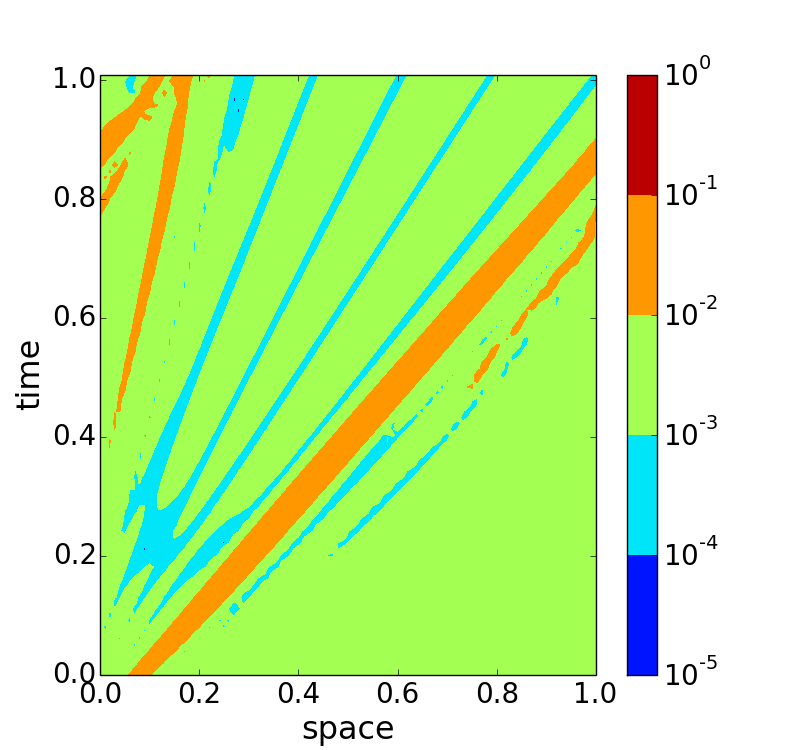}
          \end{subfigure}
    \end{subfigure}
    \begin{subfigure}[b]{.95\textwidth}
          \centering
          \begin{subfigure}[b]{0.32\textwidth}
                  \centering
                  \includegraphics[width=4.3cm]{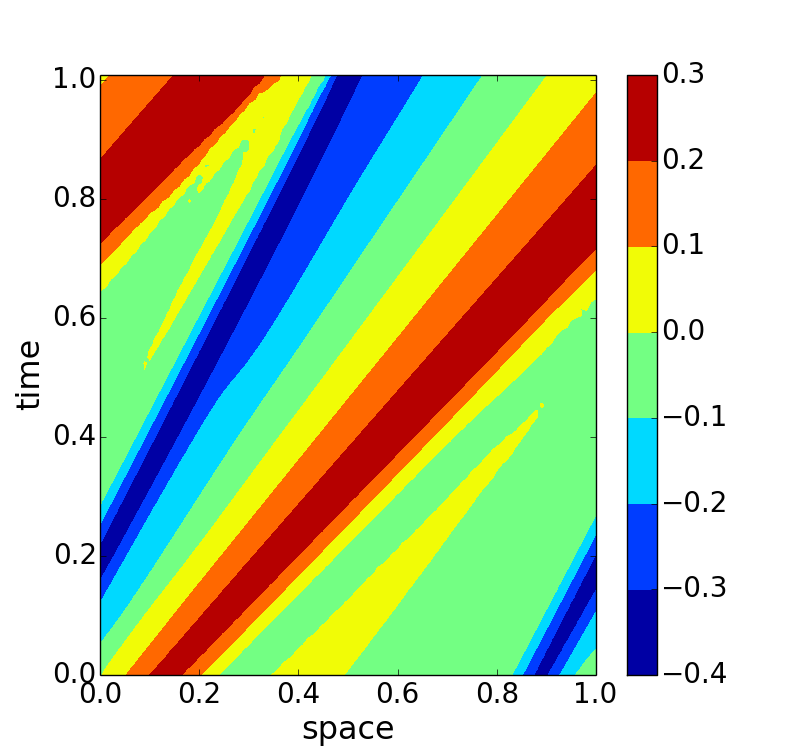}
          \end{subfigure}%
          \begin{subfigure}[b]{0.32\textwidth}
                  \centering
                  \includegraphics[width=4.3cm]{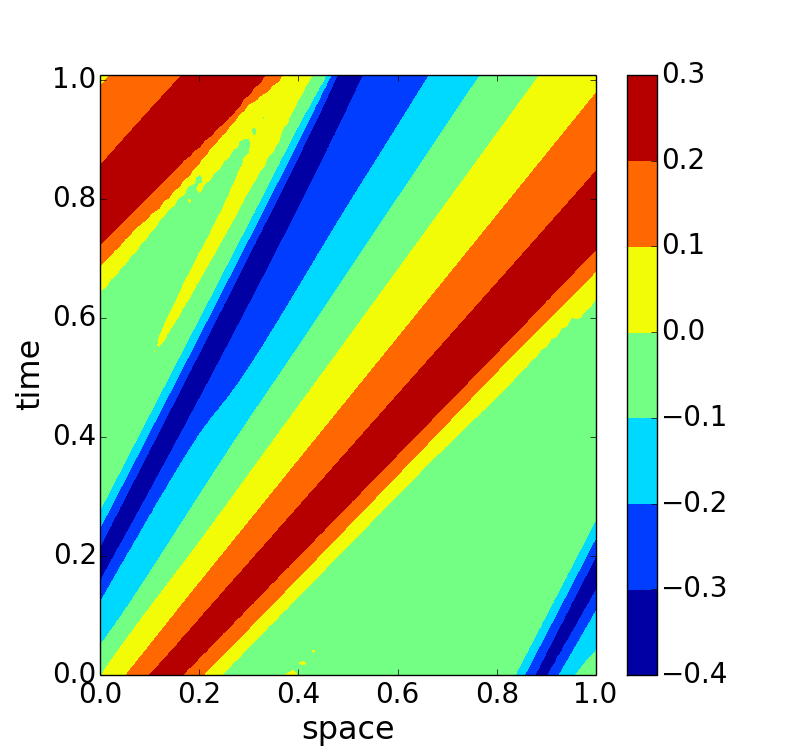}
          \end{subfigure}
          \begin{subfigure}[b]{0.32\textwidth}
                  \centering
                  \includegraphics[width=4.3cm]{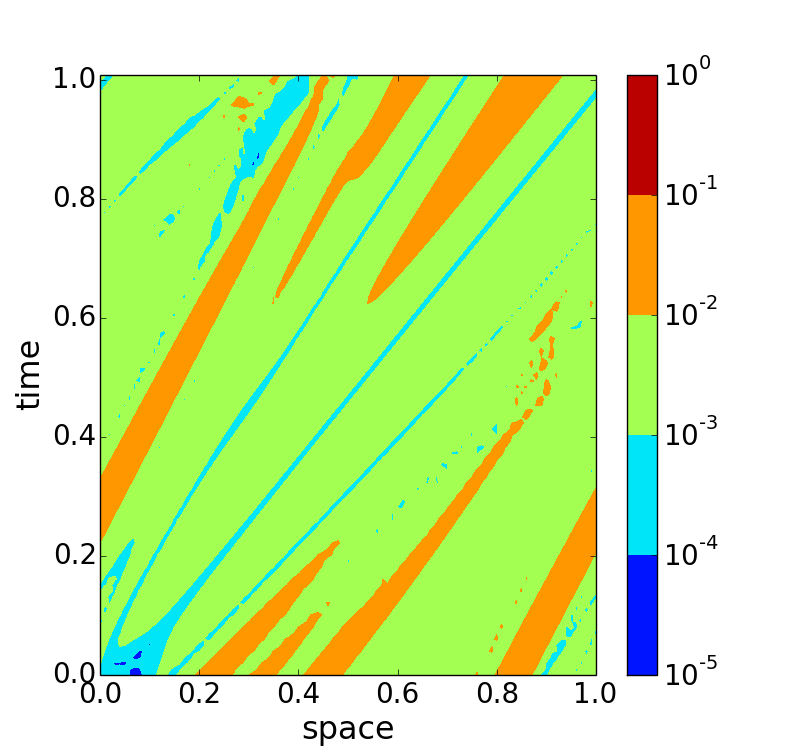}
          \end{subfigure}%
    \end{subfigure}
    \begin{subfigure}[b]{.95\textwidth}
          \centering
          \begin{subfigure}[b]{0.32\textwidth}
                  \centering
                  \includegraphics[width=4.3cm]{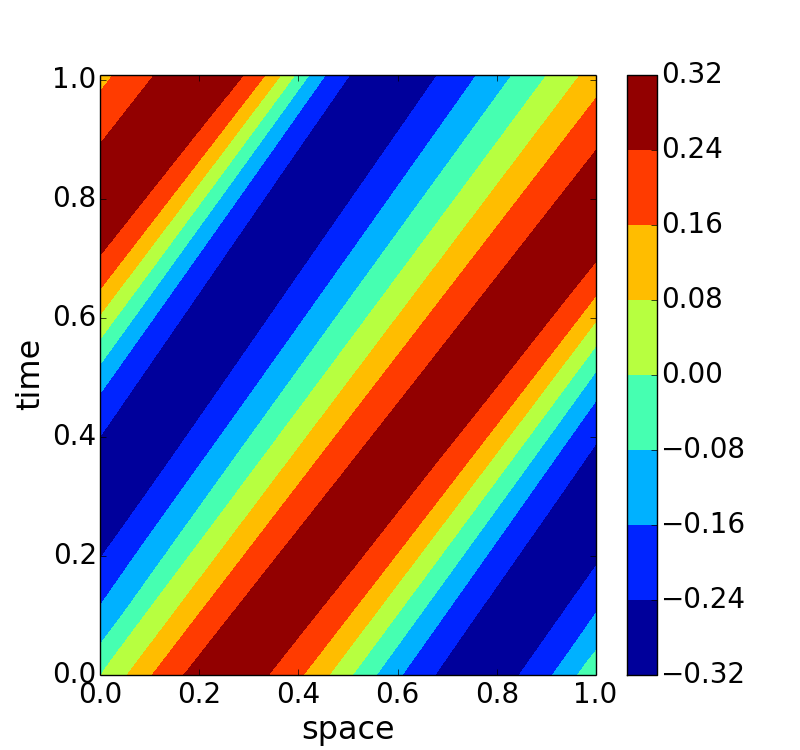}
          \end{subfigure}%
          \begin{subfigure}[b]{0.32\textwidth}
                  \centering
                  \includegraphics[width=4.3cm]{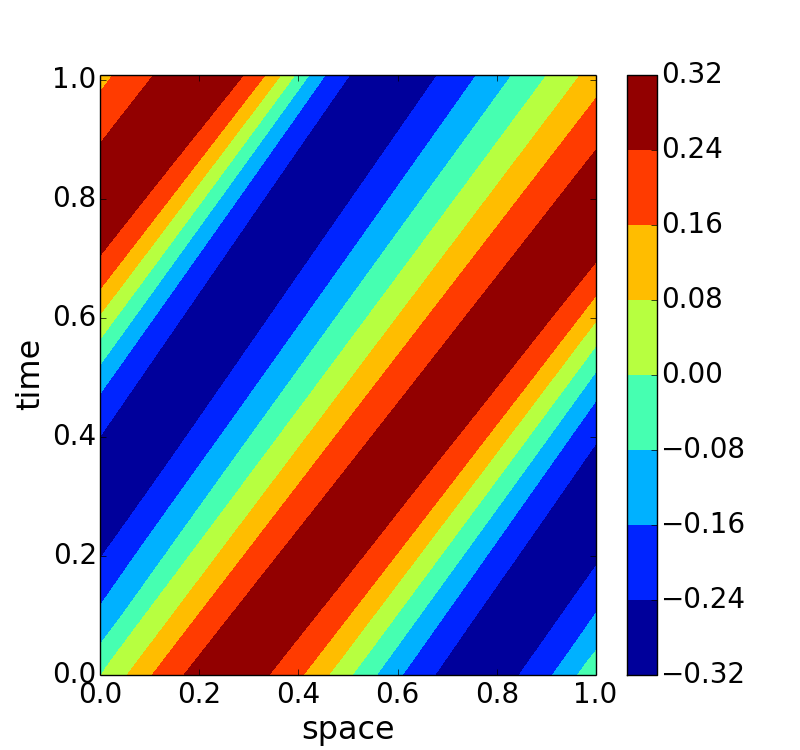}
          \end{subfigure}
          \begin{subfigure}[b]{0.32\textwidth}
                  \centering
                  \includegraphics[width=4.3cm]{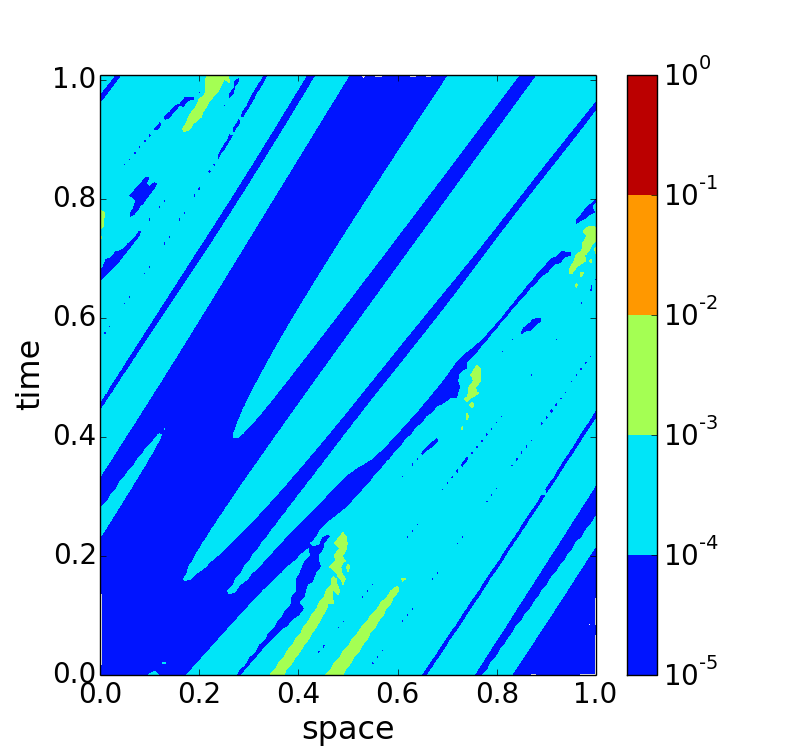}
          \end{subfigure}
    \end{subfigure}

    \caption{
             \label{fig:adj compare}
             The left column shows $\frac{dJ}{dc}$ which is evaluated by the gray-box model. 
             The middle column shows $\frac{d\tilde{J}}{dc}$ which is evaluated by the trained 
             twin model. The right column shows $\left|\frac{dJ}{dc} - \frac{d\tilde{J}}{dc}\right|$.
             We observe that the gradient is more accurate when the excited domain is smaller.}
\end{figure}

Using the twin model, we can estimate the objective's gradient by applying the
adjoint method to the twin model, i.e. we approximate
$\frac{\partial J}{\partial u}$ with $\frac{\partial \tilde{J}}{\partial \tilde{u}}$.
Suppose the gray-box model solves
\begin{equation}
    \frac{\partial u}{\partial t} + \frac{\partial}{\partial x}\,
    \left(F(u)\right) = c
\end{equation}
for $c=0$, with $F(u)$ given by Eqn.\eqref{BL flux} with $A=2$. We have trained a twin model
Eqn.\eqref{twin model 3} using the space-time solution. We are interested in
the approximation quality of the twin model's gradient at $c=0$.
$c$ is space-time dependent, therefore $\frac{dJ}{dc}$ is space-time dependent too.
We compare $\frac{dJ}{dc}$ with $\frac{d\tilde{J}}{dc}$ in Fig.\ref{fig:adj compare}

\begin{figure}
    \centering
    $\frac{dJ}{dc}$ \hspace{4.3cm} $\frac{d\tilde{J}}{dc}$  \hspace{4cm} $\left|\frac{dJ}{dc} - \frac{d\tilde{J}}{dc}\right|$\\
    \begin{subfigure}[b]{.95\textwidth}
          \centering
          \begin{subfigure}[b]{0.32\textwidth}
                  \centering
                  \includegraphics[width=4.3cm]{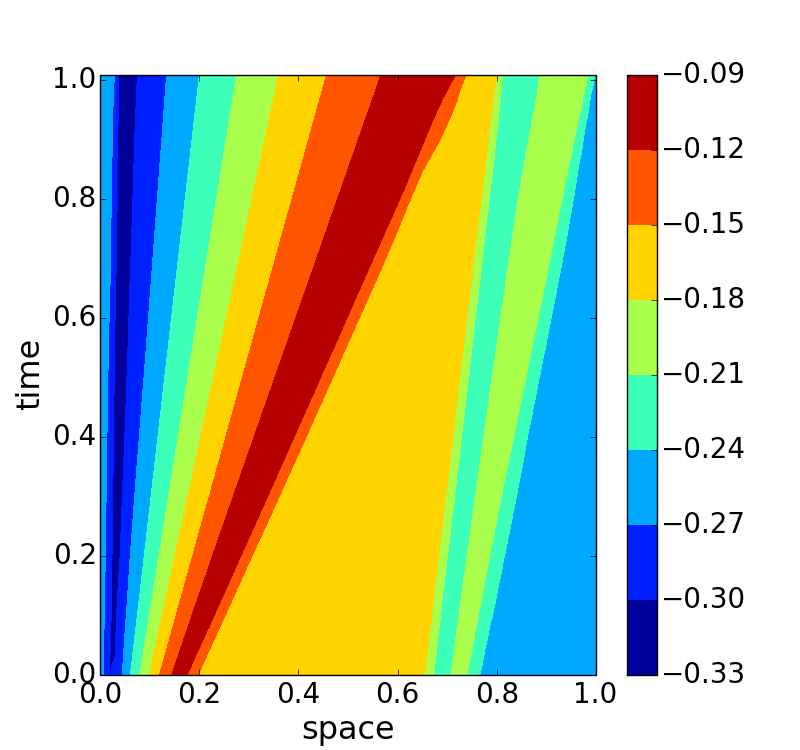}
          \end{subfigure}%
          \begin{subfigure}[b]{0.32\textwidth}
                  \centering
                  \includegraphics[width=4.3cm]{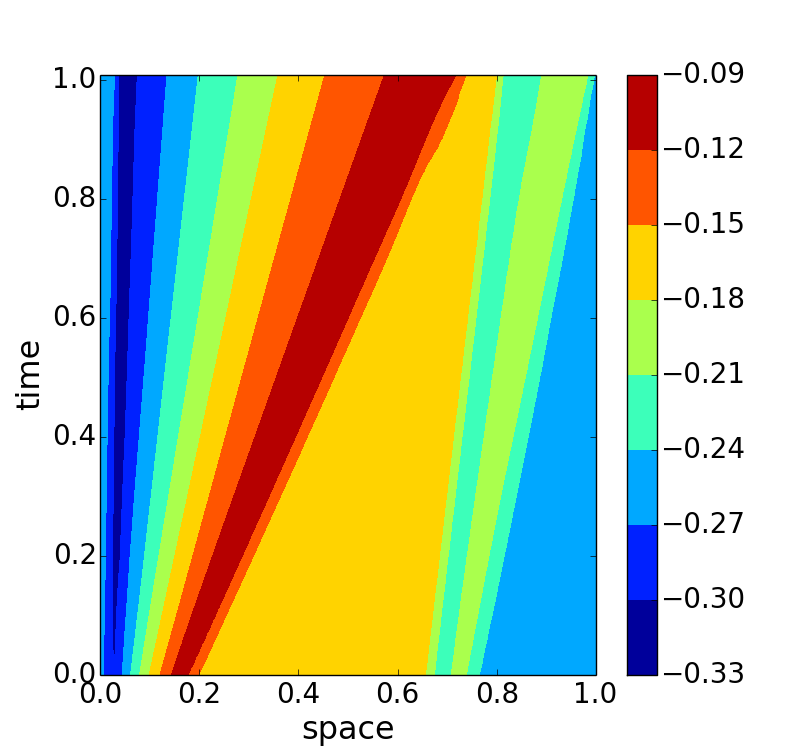}
          \end{subfigure}
          \begin{subfigure}[b]{0.32\textwidth}
                  \centering
                  \includegraphics[width=4.3cm]{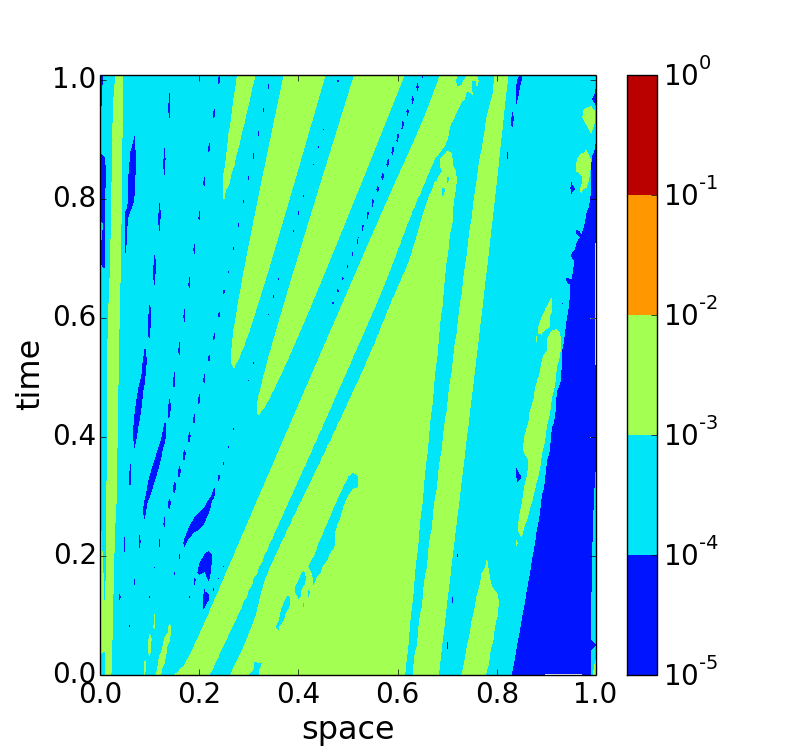}
          \end{subfigure}%
    \end{subfigure}
    \begin{subfigure}[b]{.95\textwidth}
          \centering
          \begin{subfigure}[b]{0.32\textwidth}
                  \centering
                  \includegraphics[width=4.3cm]{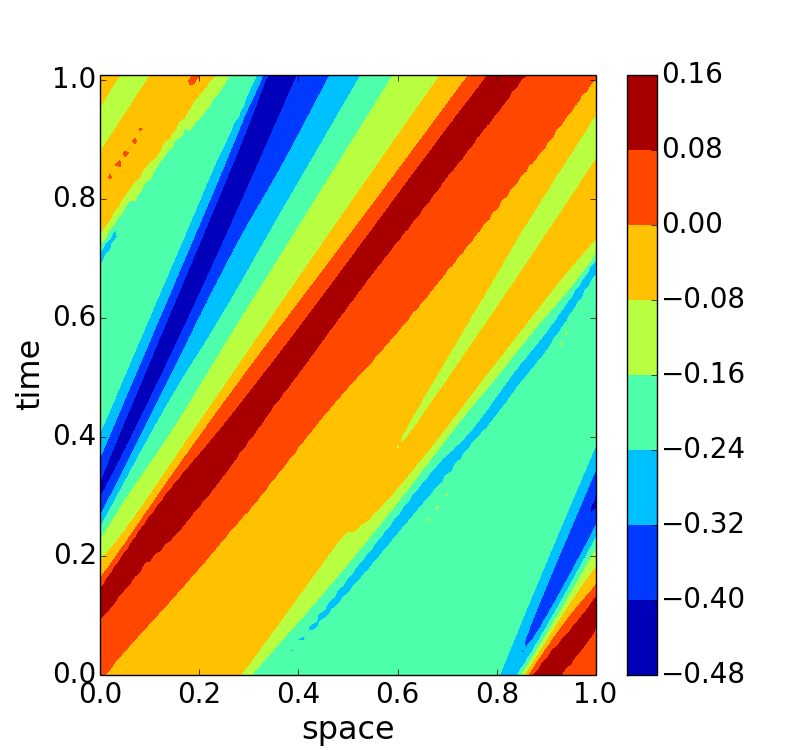}
          \end{subfigure}%
          \begin{subfigure}[b]{0.32\textwidth}
                  \centering
                  \includegraphics[width=4.3cm]{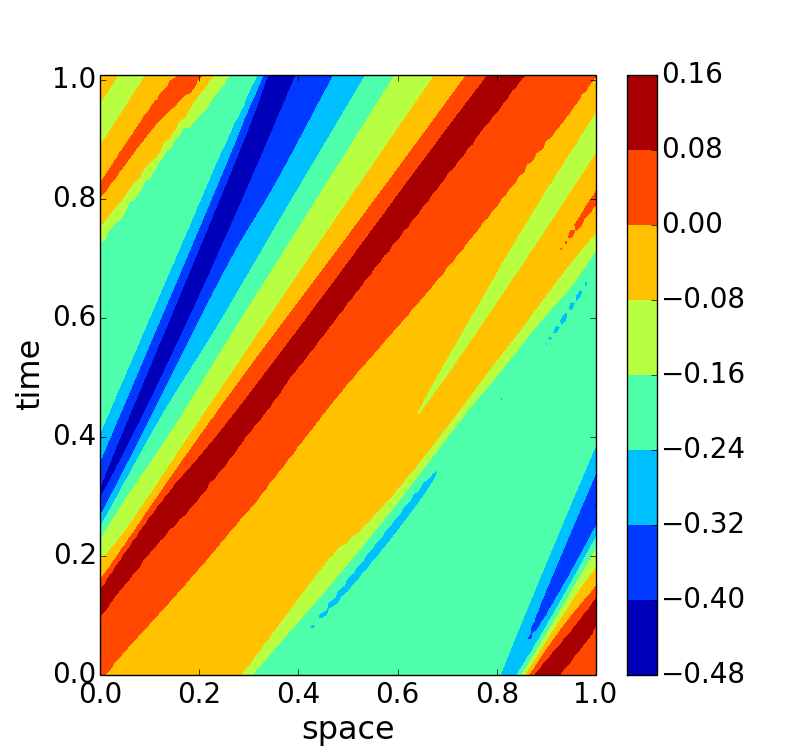}
          \end{subfigure}
          \begin{subfigure}[b]{0.32\textwidth}
                  \centering
                  \includegraphics[width=4.3cm]{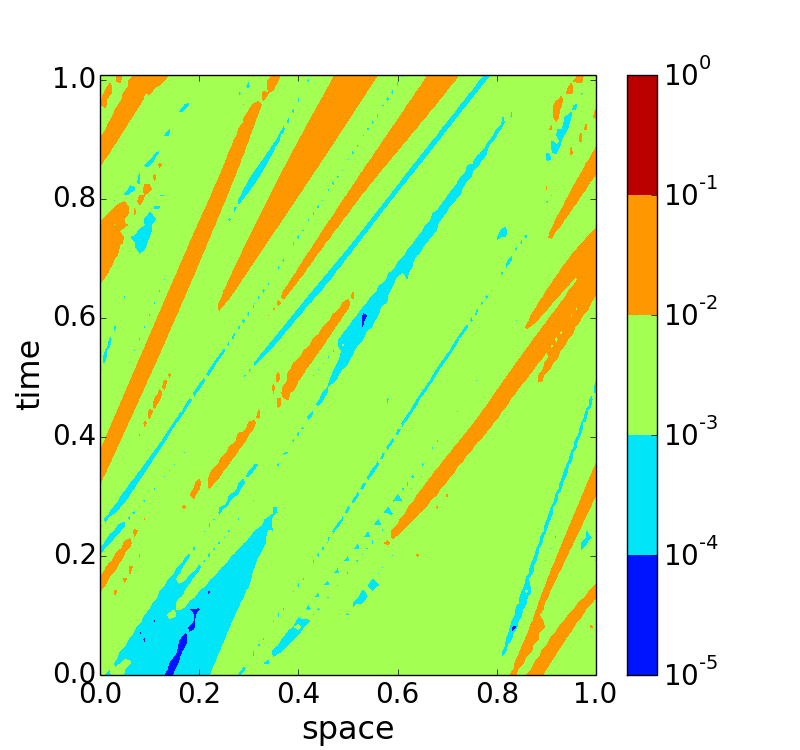}
          \end{subfigure}
    \end{subfigure}
    \begin{subfigure}[b]{.95\textwidth}
          \centering
          \begin{subfigure}[b]{0.32\textwidth}
                  \centering
                  \includegraphics[width=4.3cm]{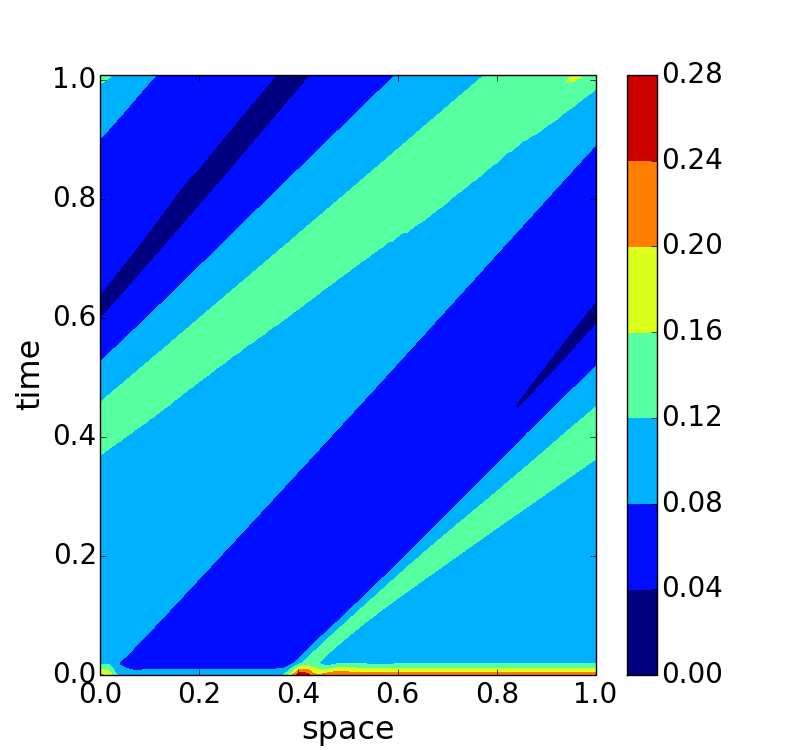}
          \end{subfigure}%
          \begin{subfigure}[b]{0.32\textwidth}
                  \centering
                  \includegraphics[width=4.3cm]{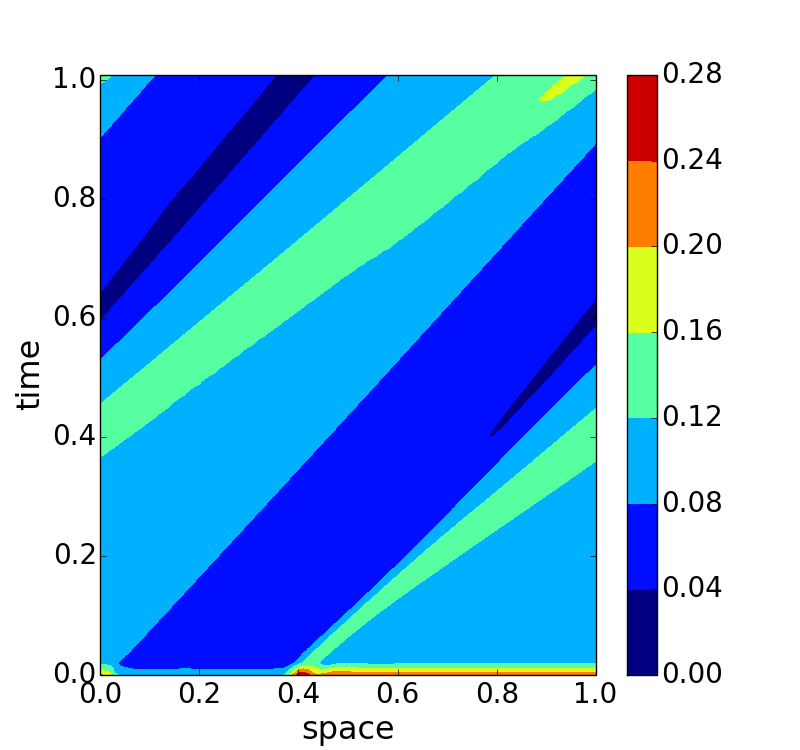}
          \end{subfigure}
          \begin{subfigure}[b]{0.32\textwidth}
                  \centering
                  \includegraphics[width=4.3cm]{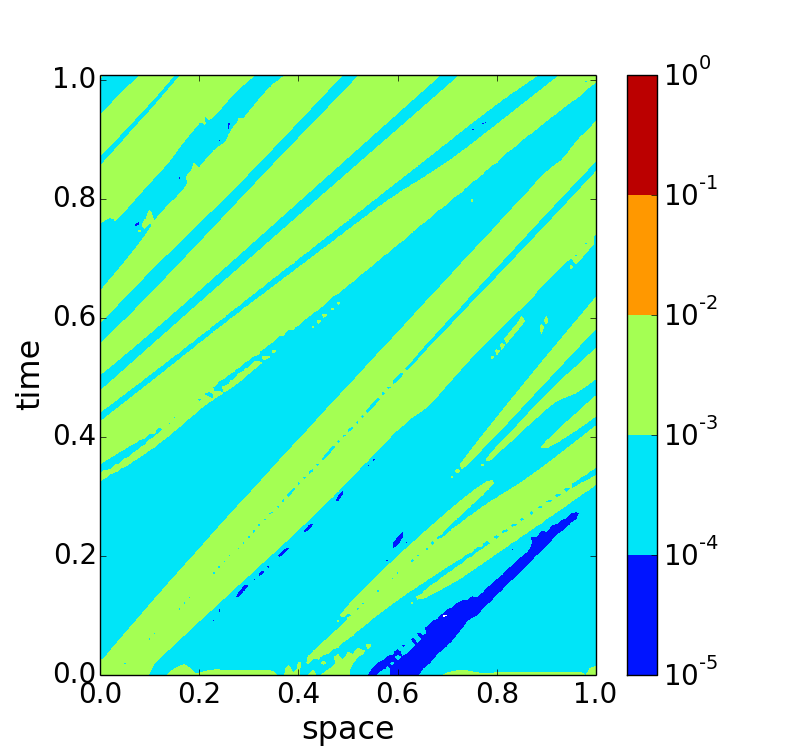}
          \end{subfigure}%
    \end{subfigure}
    \begin{subfigure}[b]{.95\textwidth}
          \centering
          \begin{subfigure}[b]{0.32\textwidth}
                  \centering
                  \includegraphics[width=4.3cm]{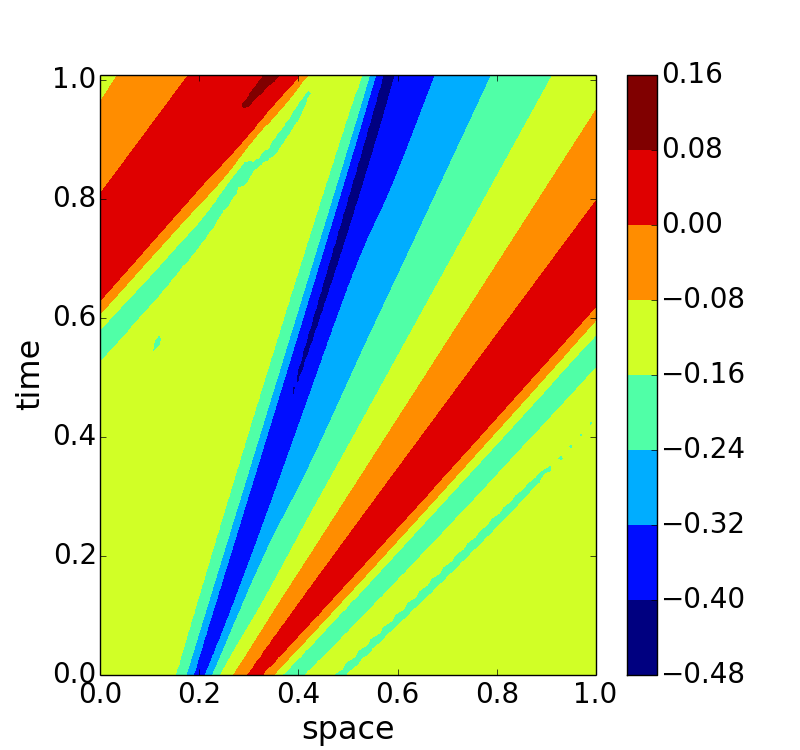}
          \end{subfigure}%
          \begin{subfigure}[b]{0.32\textwidth}
                  \centering
                  \includegraphics[width=4.3cm]{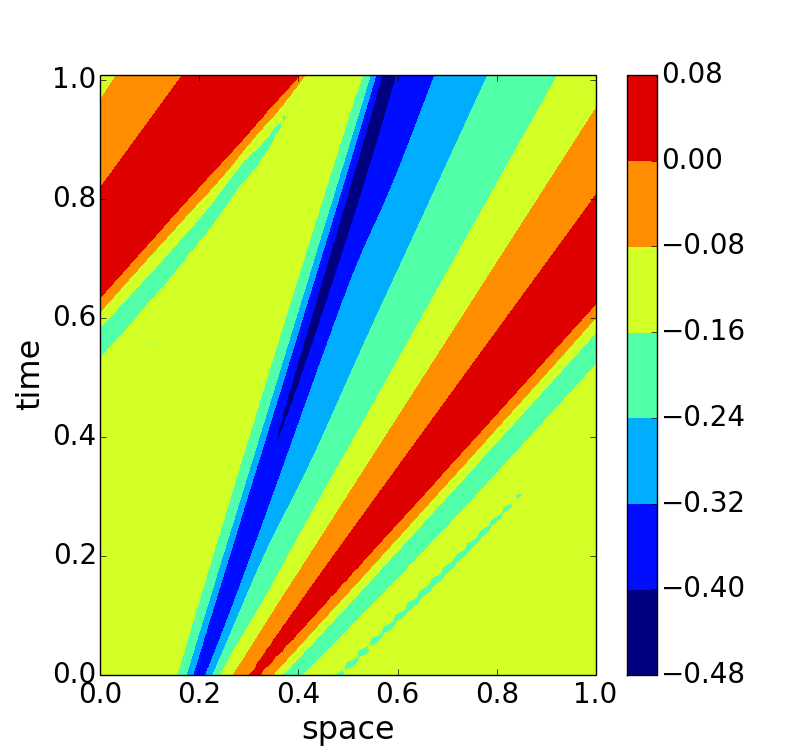}
          \end{subfigure}
          \begin{subfigure}[b]{0.32\textwidth}
                  \centering
                  \includegraphics[width=4.3cm]{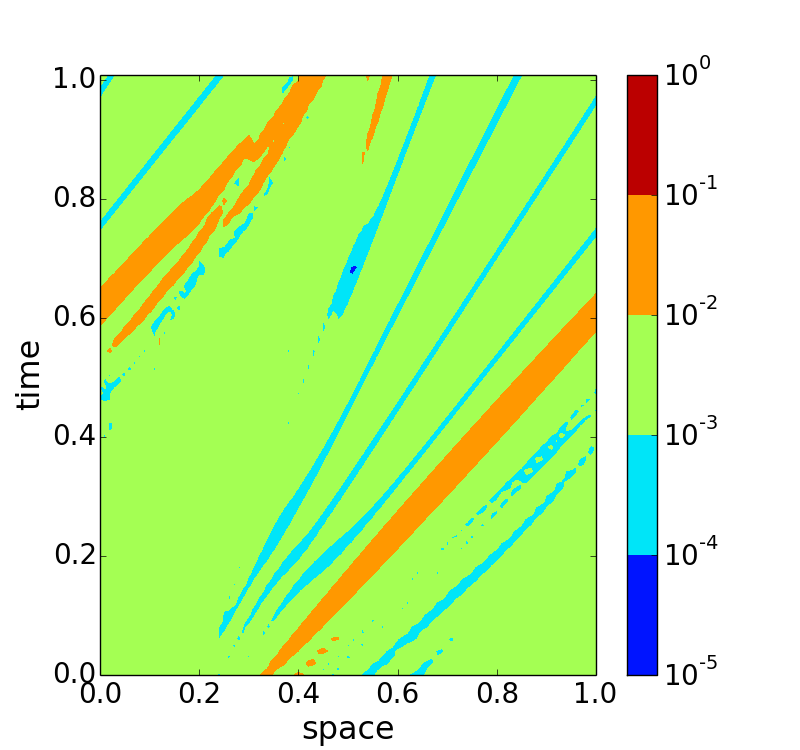}
          \end{subfigure}
    \end{subfigure}
    \begin{subfigure}[b]{.95\textwidth}
          \centering
          \begin{subfigure}[b]{0.32\textwidth}
                  \centering
                  \includegraphics[width=4.3cm]{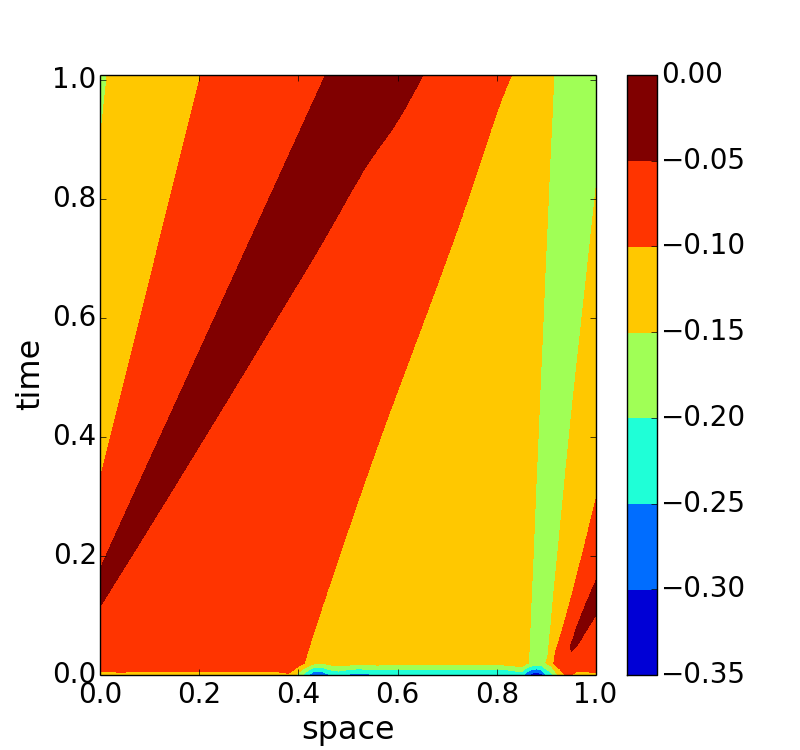}
          \end{subfigure}%
          \begin{subfigure}[b]{0.32\textwidth}
                  \centering
                  \includegraphics[width=4.3cm]{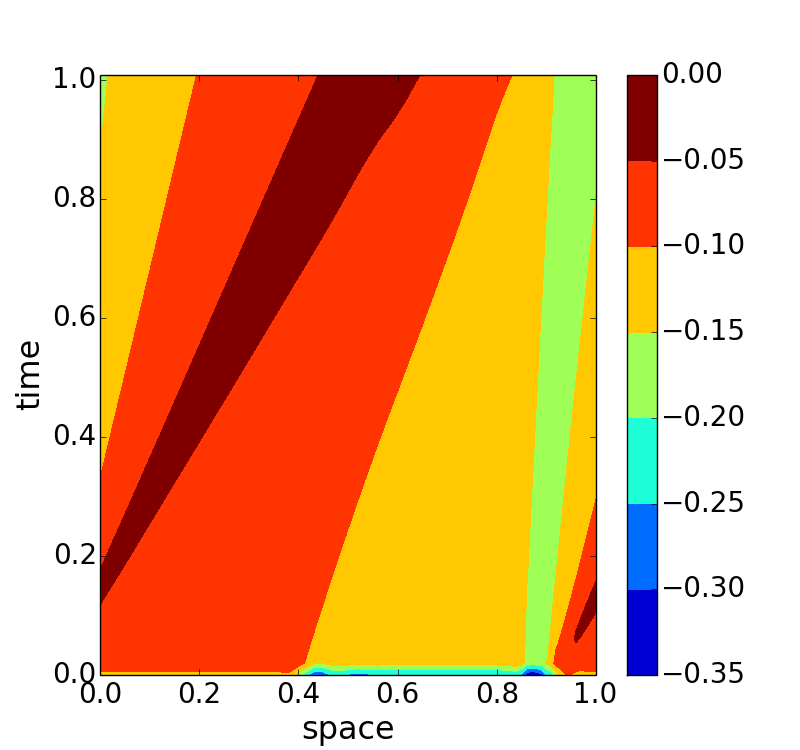}
          \end{subfigure}
          \begin{subfigure}[b]{0.32\textwidth}
                  \centering
                  \includegraphics[width=4.3cm]{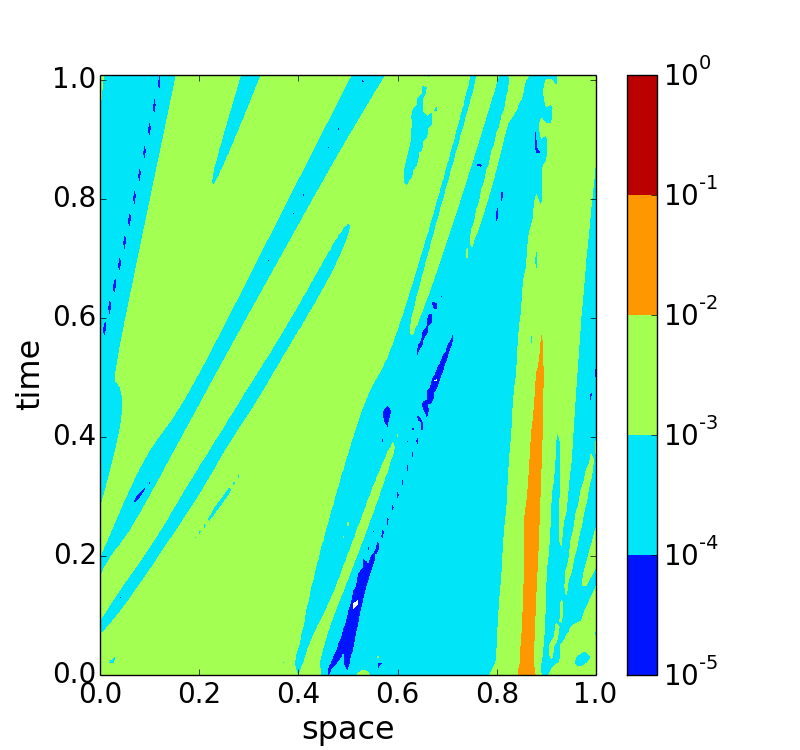}
          \end{subfigure}%
    \end{subfigure}
\end{figure}

The result is encouraging, as the gradient computed by the twin model provides
a good approximation
of the gradient of the primal model. We reiterate that the good approximation quality benefits from 
the matching of the space-time solution.

\clearpage
\subsection{Gradient estimation of Navier-Stokes flows}
In the second numerical test case, we consider
a compressible internal flow in a 2-D return bend channel.
The flow is driven by the pressure difference between the inlet and the outlet.
The flow is governed by Navier-Stokes equations, Eqn.\eqref{NSeqn}.
Navier-Stokes equations require an additional state equation, Eqn.\eqref{state equation}, for closure.
Many models of the state equations have been developed, including the ideal gas equation, the
van der Waals equation, and the Redlich-Kwong equation \cite{state eqns}.\\

\begin{figure}\begin{center}
    \includegraphics[width=10cm]{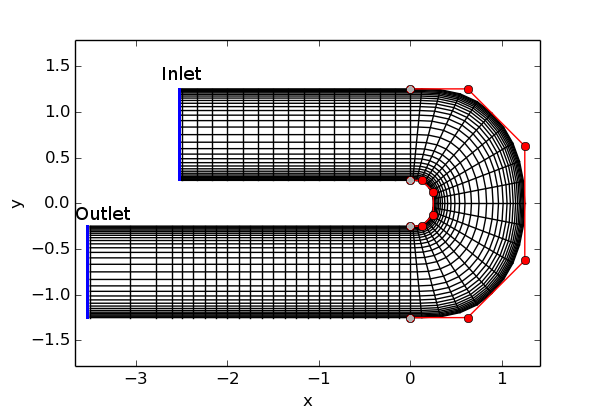}
    \caption{The return bend geometry and the mesh for simulation. The return bend is bounded by 
    no-slip walls. The values of pressure
    are fixed at the inlet and at the outlet. 
    The inner and outer boundaries of the bend are generated by the control points
    using quadratic B-spline.
    We estimate the gradient of the steady state mass flux
    to the red control points' coordinates.}
    \label{NS mesh}
\end{center}\end{figure}

The inner and outer boundaries at the bending section 
are each generated by 6 control points using quadratic B-spline. 
The control points are shown by the red and gray dots in Fig.\ref{NS mesh}.
The gray control points are fixed 
on the straight sections.
The spanwise grid is generated by geometric grading. The streamwise
grid at the straight and the bending section are each generated
by uniform grading, except at the sponge region.
The pressure at the outlet is set to be a constant $p_{out}$ while
the total pressure at the inlet is set to be a constant $p_{t,in}$.
Let $\rho_\infty$ be the steady state density, and
$\boldsymbol{u}_\infty = (u_\infty, v_\infty)$ be the steady state Cartesian velocity.
The steady state mass flux is
\begin{equation}
    J = - \int_{\textrm{outlet}} \rho_\infty u_\infty \big|_{\textrm{outlet}} \; dy=
    \int_{\textrm{inlet}} \rho_\infty u_\infty\big|_{\textrm{inlet}} \; dy
    \label{mass flux}
\end{equation}
We want to estimate the gradient of the steady state mass flux
to the red control points' coordinates.
\\

\begin{figure}\begin{center}
    \includegraphics[width=14cm]{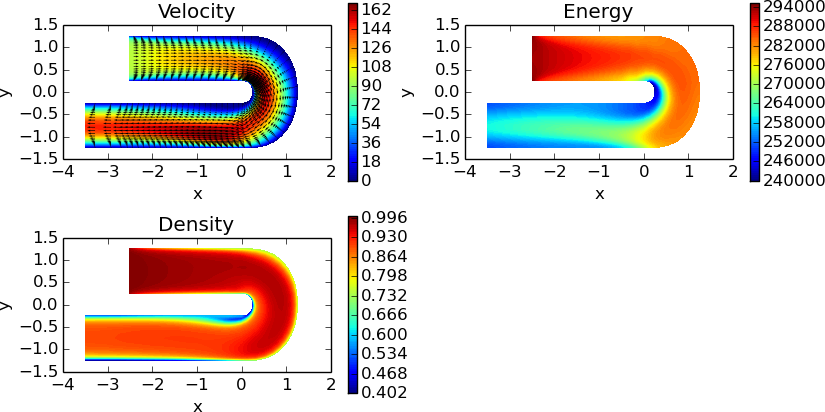}
    \caption{An example of the steady state velocity, energy, and density provided 
    by the gray-box simulation. In the velocity subplot, the 
    magnitude of the velocity is overlayed with the velocity vectors.}
    \label{gray sol}
\end{center}\end{figure}

When the state equation of the fluid is unknown, the 
adjoint method cannot be applied directly
to estimate the gradient. We use the proposed twin model to infer the state
equation from the steady state solution of the gray-box simulation.
Assume that the gray-box simulation provides $\rho_\infty$, $\boldsymbol{u}_\infty$,
and $E_\infty$. 
Fig.\ref{gray sol} shows an example of the solution.

We parameterize the unknown state equation by
\begin{equation}
    p(\rho, U) = \sum_{i=1}^{N_\rho} \sum_{j=1}^{N_U} \alpha_{ij} R_i(\rho) S_j(U) 
    + p_0
    \label{parameterization}
\end{equation}
where
\begin{equation}
    R_i(\rho) = \exp\left(\frac{-(\rho-\rho_i)^2}{\sigma_\rho}\right)\,,\quad i=1,\cdots,N_\rho 
\end{equation}
\begin{equation}
    S_j(U) = \frac{1}{2}\left(\tanh\left( \frac{U-U_j}{\sigma_U}\right) +1\right)
    \quad j = 1,\cdots,N_U\,.
\end{equation}
$R_i$ are radial basis functions, and $S_j$ are sigmoid functions.
Let the density of the gray-box solution be in the range of $[\rho_{\min}, \rho_{\max}]$, and
let the internal energy of the gray-box solution be in the range of $[U_{\min}, U_{\max}]$.
We set $\rho_i$, $i=1,\cdots,N_\rho$ to be equally spaced in $[\rho_{\min}, \rho_{\max}]$, and set
$U_j$, $j=1,\cdots,N_U$ to be equally spaced in $[U_{\min}, U_{\max}]$.
We constrain $\alpha$ values to be positive to respect the fact that pressure
monotonically increases with the internal energy.
$p_0$ is a scalar. We set $\sigma_\rho=\frac{\rho_{\max}-\rho_{\min}}{N_\rho}$, 
$\sigma_U = \frac{U_{\max}-U_{\min}}{N_U}$.\\

We define the solution mismatch, Eqn.\eqref{minimizer twin model discrete steady}, as
\begin{equation}
    \mathcal{M} = w_\rho \|\tilde{\rho}_{\infty} - \rho_{\infty}\|^2
                + w_u
                \|\tilde{u}_{\infty}- u_{\infty}\|^2
                + w_v
                \|\tilde{v}_{\infty}- v_{\infty}\|^2
                + w_E
                \|\tilde{E}_{\infty} - E_\infty\|^2
    \label{NS mismatch}
\end{equation}
where $w_\rho$, $w_u$, $w_v$, and $w_E$ are positive weight constants.
$\|\cdot\|$ is the $L_2$ norm.
To infer the state equation we solve the optimization problem:
\begin{equation}
    \min_{\alpha, p_0} \left\{\mathcal{M} +\lambda 
    \sum_{i=1}^{N_\rho}\sum_{j=1}^{N_U} \big|\alpha_{ij}\big|
    \right\}
    \label{NS optimize}
\end{equation}

We need to choose suitable weights $w_\rho$, $w_u$, $w_v$, and $w_E$ in
Eqn.\eqref{NS optimize}. To select these weights, 
we first set $\alpha$ and $p_0$ values to several randomly guessed values. 
Using these guessed state equation, we 
obtain $\|\tilde{\rho}_{\infty} - \rho_{\infty}\|^2$,
$\|\tilde{u}_{\infty}- u_{\infty}\|^2 $,
$ \|\tilde{v}_{\infty}- v_{\infty}\|^2$,
and $\|\tilde{E}_{\infty} - E_\infty\|^2$.
The weights are chosen to be
\begin{equation}\begin{split}
    w_\rho &= \frac{1}{\left<\|\tilde{\rho}_{\infty} - \rho_{\infty}\|^2\right>}\\
    w_u &= \frac{1}{\left<\|\tilde{u}_{\infty} - 
    {u}_{\infty}\|^2\right>}\\
    w_v &= \frac{1}{\left<\|\tilde{v}_{\infty} - 
    v_{\infty}\|^2\right>}\\
    w_E &= \frac{1}{\left<\|\tilde{E}_{\infty} - 
    E_{\infty}\|^2\right>}
\end{split}
\label{NS weights}
\end{equation}
where $\left<\cdot\right>$ denotes the sample average of the randomly-guessed state equations.
In this way, $\tilde{\rho}_\infty-\rho_\infty$, 
$\tilde{u}_{\infty}-u_\infty$, $\tilde{v}_{\infty}-v_\infty$,
and $\tilde{E}_\infty-E_\infty$ will, on average, contribute
similarly to the solution mismatch for the randomly-guessed state equations.\\

We tested three example state equations in the graybox simulator:
the ideal gas equation, the van der Waals equation, and the Redlich-Kwong equation:
\begin{equation}\begin{split}
    p_{ig} &= (\gamma-1) U\\
    p_{vdw} &= \frac{(\gamma-1)U}{1-b_{vdw}\rho} - a_{vdw}\rho^2\\
    p_{rk} &= \frac{(\gamma-1)U}{1-b_{rk}\rho} - 
    \frac{a_{rk}\rho^{5/2}}{((\gamma-1)U)^{1/2}(1+b_{rk}\rho)}
\end{split}\label{NS state equations}
\end{equation}
where $a_{vdw}$, $b_{vdw}$, $a_{rk}$, $b_{rk}$ are constants. In the following testcases, we choose
$a_{vdw}=10^4$, $b_{vdw}=0.1$, $a_{rk}=10^7$, $b_{rk}=0.1$.\\

\begin{figure}
    \centering
    Redlich-Kwong gas\\
    \includegraphics[width=14cm]{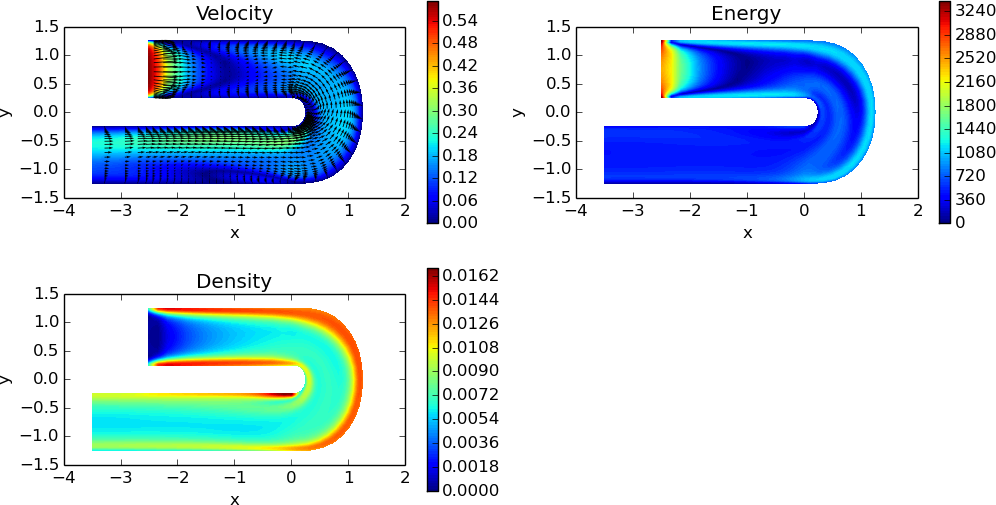}\\
    \caption{
    \label{fig:NS sol err}
    The solution mismatch between the twin model's solution and the graybox solution for the
    ideal gas equation, the van der Waals equation, and the Redlich-Kwong equation. 
    The first group of images is for the ideal gas, the second group of images is for the van der Waals
    gas, and the third group of images is for the Redlich-Kwong gas. For each gas, we show 
    the solution mismatch of velocity, energy, and density. For all
    three test cases, the relative error of the velocity, energy, and density is small.}
\end{figure}

By solving Eqn \eqref{NS optimize}, we obtain the solution mismatch for the state equations.
Fig.\ref{fig:NS sol err} shows the solution mismatch,
 $\left|\tilde{\rho} -\rho\right|$, $\left|\tilde{\boldsymbol{u}}- \boldsymbol{u}\right|$, 
and $\left|\tilde{E}-E\right|$
from the inferred twin model and the graybox model.
\begin{figure}
    \centering
    Ideal gas\\
    \includegraphics[width=14cm]{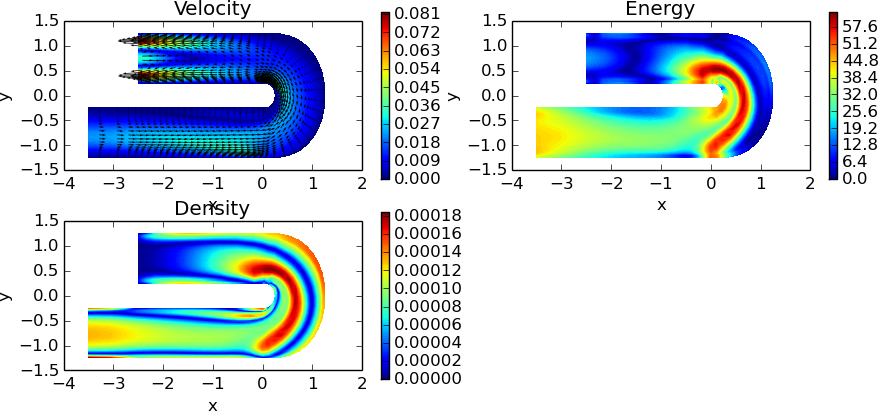}\\
    van der Waals gas\\
    \includegraphics[width=14cm]{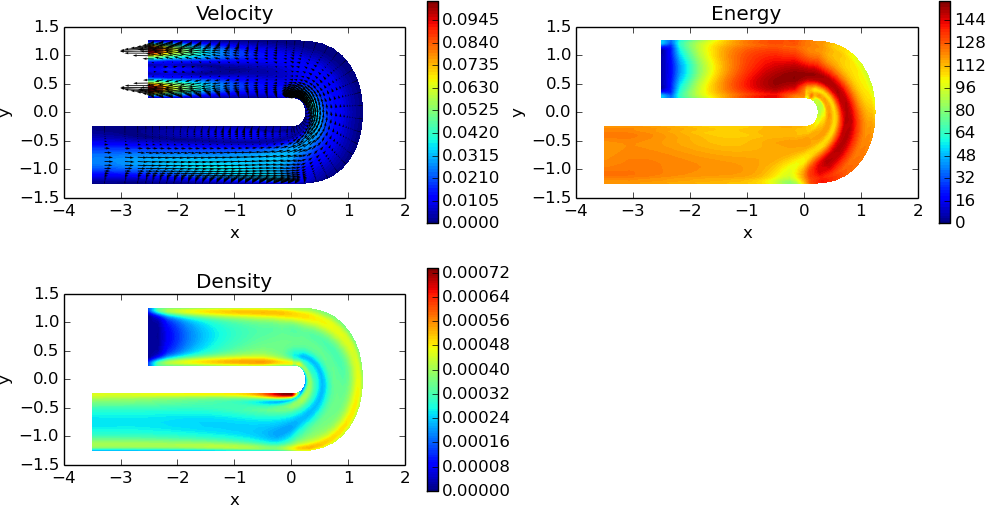}\\
\end{figure}

\begin{figure}
\begin{center}
    Redlich-Kwong gas\\
    \includegraphics[width=14cm]{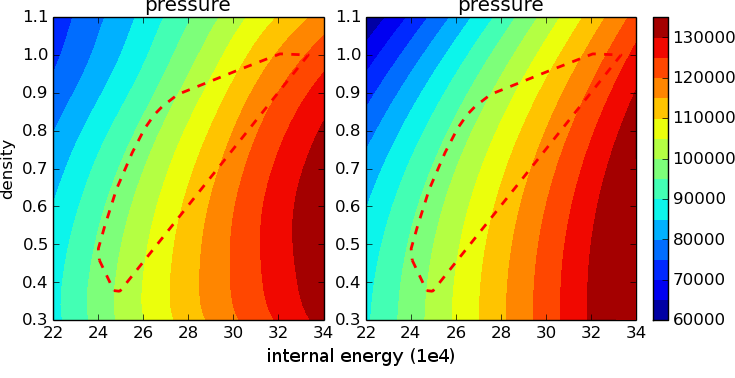}
    \caption{
    \label{fig:NS state eqn}
    The inferred state equation and the gray-box state equation for the three gases. 
    The left column shows the inferred state equation, and the right column shows the 
    gray-box state equation.
    In all the three test cases, the inferred state equation approximates the gray-box state
    equation accurately inside $H(\rho, U)$ indicated by the dashed line.}
\end{center}
\end{figure}

Let the $(\rho, U)$ be the gray-box steady state solution's 
density and internal energy at all the spatial
gridpoint, and let $H(\rho, U)$ be its convex hull. We expect the the estimated state equation
to be more accurate inside $H(\rho, U)$ than outside $H(\rho, U)$.
The inferred state equation, $p = p(\rho, U)$, is shown in Fig.\ref{fig:NS state eqn}.
\clearpage
\begin{figure}
\begin{center}
    Ideal gas\\
    \includegraphics[width=14cm]{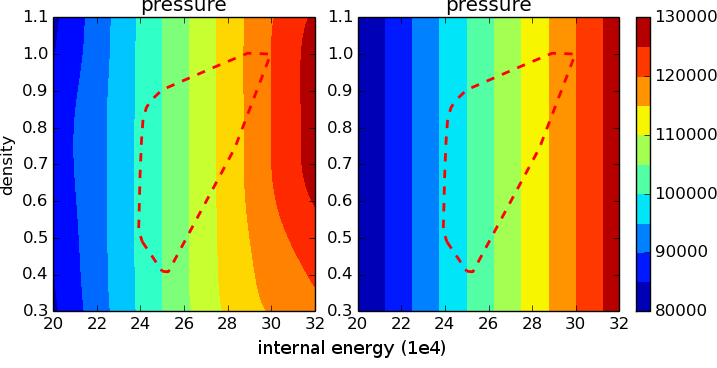}
    van der Waals gas\\
    \includegraphics[width=14cm]{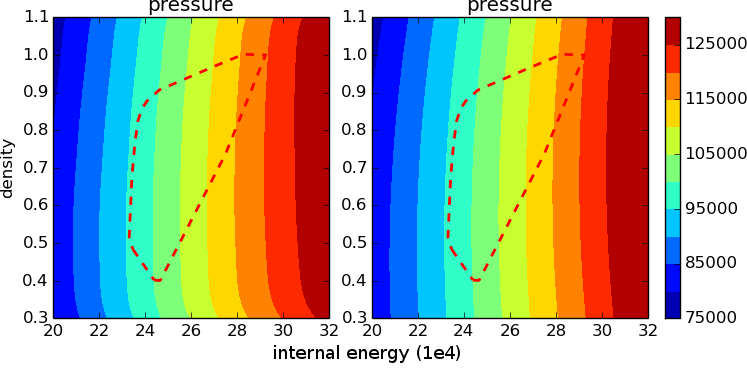}
\end{center}
\end{figure}

\begin{figure}\begin{center}
    \includegraphics[width=5.2cm]{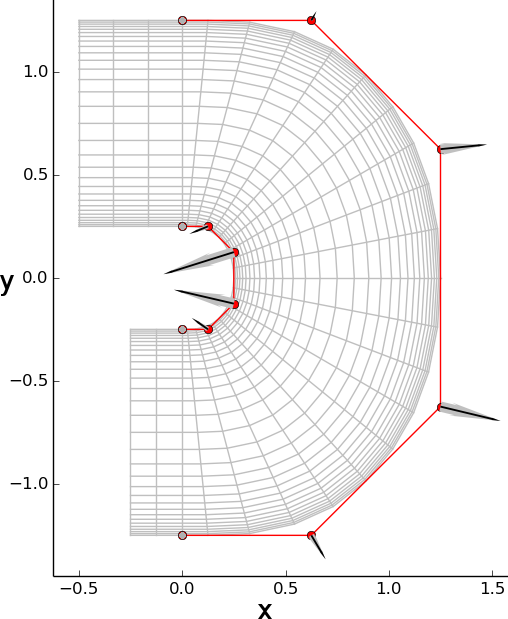}
    \includegraphics[width=5.2cm]{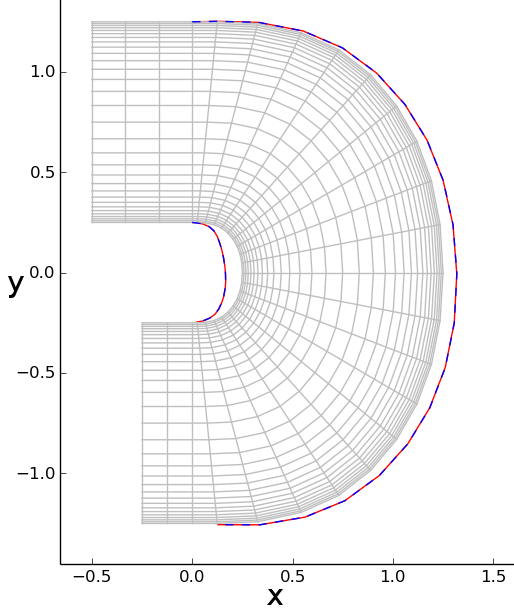}
    \caption{\label{fig:NS perturb geo}
    The left column shows the gradient of the outflux to the control points for the 
    ideal gas. 
    The wide gray arrow is the gradient evaluated by the gray-box model, while
    the thin black arrow is the gradient evaluated by the twin model.
    The right column shows a perturbed boundary according to the gradient. 
    The blue dashed line is computed by the gray-box model's gradient, while the
    red dashed line is computed by the twin model's gradient.}
\end{center}\end{figure}

\clearpage
Using the inferred state equation, we are able to compute the
gradient of the mass flux to the countrol points at the bending section.
For example, the gradient and the perturbed boundary for the ideal gas are shown in 
Fig.\ref{fig:NS perturb geo}.
For all the three gases, the difference between the gray-box gradient and the
twin model gradient is hardly visible.\\ 

\begin{figure}\centering
    \includegraphics[width=12cm]{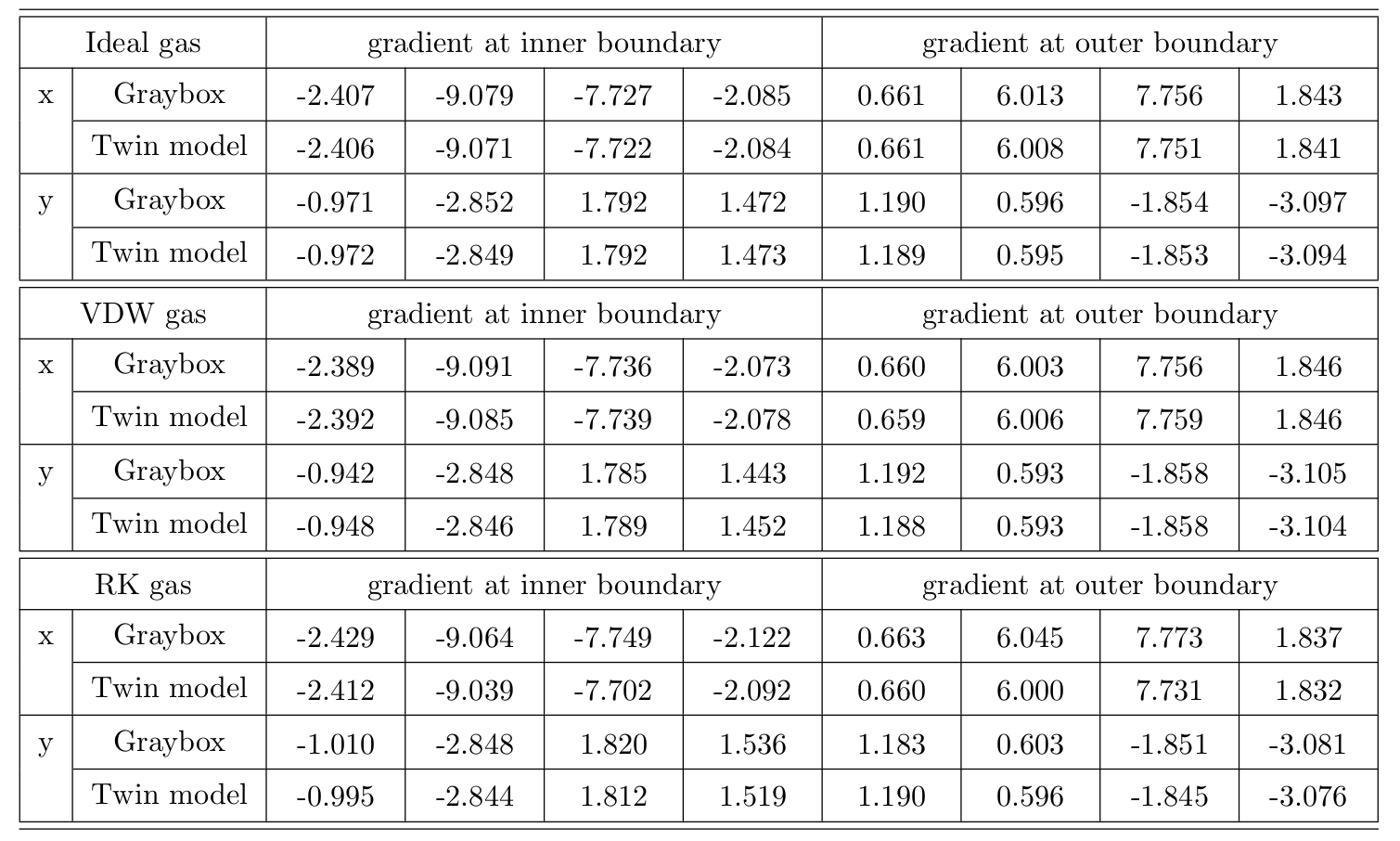}
    \caption{The gradient of the mass flux to the control points' coordinates}
    \label{tab: idea gas gradient}
\end{figure}

We summarize the estimated gradient computed by the twin model
in Fig.\ref{tab: idea gas gradient}, which is compared
with the gradient computed by the gray-box model. For each gas,
we compare the x-component and the y-component of the gradients.
Twin model demonstrates to estimate the gradients accurately.

\clearpage

%
%
%

\section{Conclusion}
We propose a method to estimate the objective's gradient when the simulator is a gray-box conservation
law simulator and does not
implement the adjoint method. 
The proposed method uses the space-time or spatial solution of the gray-box simulation to
infer a twin model. 
There are several benefits to use the space-time or spatial solution. Firstly, in many conservation law
simulations, flow quantities have a small domain of dependence. Secondly, the space-time or spatial
solution from a single simulation provides a large number of samples for the inference. 
Thirdly, in many high-dimensional design problems, the design variables are space-time or spatially distributed, 
so the inference's input dimension does not scale up with the design dimension.
The twin model method enables adjoint computation. We use the gradient computed 
by the twin model to estimate the gradient of the gray-box simulation.\\

The twin model method is demonstrated on a 1-D porous media flow problem and a
2-D Navier-Stokes problem. In the 1-D problem  
the flux function is unknown. We are able to infer the flux function in the excited domain. 
Using the inferred twin model, we estimate the gradient of the objective to the 
space-time dependent control.
In the 2-D problem, the state equation is unknown. We are able to infer
the state equation using the steady state solution of the gray-box model.
Using the inferred state equation, we estimate the gradient of the mass flux
to the coordinates of the control points.
Our research shows that gradient can be efficiently estimated by adjoint method even if 
the simulator is gray-box.\\

The twin model enables adjoint computation for gray-box simulations. 
In the future, we plan to apply the twin model to high-dimensional optimization problems.

\clearpage
\bibliographystyle{ctr}
\bibliography{ctr_summer_lesopt}

\end{document}